\documentclass[11pt,reqno]{amsart}

\usepackage{geometry}
\usepackage{amsmath,amsthm,amssymb,amsfonts,dsfont}
\usepackage{mathtools}
\usepackage{float}
\usepackage[utf8]{inputenc} 
\usepackage[T1]{fontenc}    
\usepackage[hidelinks]{hyperref}
\usepackage[capitalise]{cleveref}
\usepackage{pgfplots}
\usepackage{comment}
\usepackage{subcaption}
\usepackage{graphicx}
\usepackage{wrapfig}
\usepackage{bm}
\usepackage{url}           
\usepackage{booktabs}    
\usepackage{nicefrac}    
\usepackage{microtype}    
\usepackage{bbm}
\usepackage{psfrag}
\usepackage{breakcites}
\usepackage{color,soul}
\usepackage{breakcites}   
\usepackage{bigints}
\usepackage{enumitem}
\setlist{leftmargin=1.6em}

\newtheorem{remark}{Remark}
\newtheorem{theorem}{Theorem}
\newtheorem{proposition}{Proposition}
\newtheorem{lemma}{Lemma} 

\usepackage{graphicx}
\usepackage{xcolor}
\usepackage{svg}
\usepackage{algorithm}
\usepackage{cite}
\usepackage[noend]{algpseudocode}
\hypersetup{colorlinks, linkcolor={red!80!black}, citecolor={blue!80!black}, urlcolor={ForestGreen}}
\graphicspath{{figs/}}

\makeatletter
\def\subsubsection{\@startsection{subsubsection}{3}%
  \z@{.5\linespacing\@plus.7\linespacing}{-.5em}%
  {\normalfont\bfseries}}
\makeatother

\newcommand{\sign}[0]{\mathrm{sign}}

\newcommand{\calX}[0]{\mathcal{X}}

\newcommand{\calY}[0]{\mathcal{Y}}

\newcommand{\R}[0]{\mathbb{R}}

\newcommand{\vertiii}[1]{{\left\vert\kern-0.25ex\left\vert\kern-0.25ex\left\vert #1 
    \right\vert\kern-0.25ex\right\vert\kern-0.25ex\right\vert}}

\newcommand{\sF}{\mathsf{F}}
\newcommand{\sH}{\mathsf{H}}

\newcommand{\sC}{\mathsf{C}}
\newcommand{\bA}{\mathbf{A}}
\newcommand{\bB}{\mathbf{B}}

\newcommand{\cC}{\mathcal{C}}
\newcommand{\cD}{\mathcal{D}}

\newcommand{\cF}{\mathcal{F}}

\newcommand{\cP}{\mathcal{P}}

\newcommand{\cU}{\mathcal{U}}

\newcommand{\cX}{\mathcal{X}}
\newcommand{\cY}{\mathcal{Y}}

\newcommand{\EE}{\mathbb{E}}

\newcommand{\NN}{\mathbb{N}}

\newcommand{\RR}{\mathbb{R}}

\newcommand{\vasti}{\bBigg@{3.5 }}
\newcommand{\vast}{\bBigg@{4}}
\newcommand{\Vast}{\bBigg@{5}}
\newcommand{\Vastt}{\bBigg@{7}}

\makeatother

\newcommand{\be}{\begin{equation}}
\newcommand{\ee}{\end{equation}}
\newcommand{\ba}{\begin{align}}
\newcommand{\ea}{\end{align}}
\newcommand{\baa}{\begin{align*}}
\newcommand{\eaa}{\end{align*}}

\newcommand{\KL}{\mathsf{D}_{\mathsf{KL}}}
\newcommand{\OT}{\mathsf{OT}}
\newcommand{\IOT}{\mathsf{IOT}}

\DeclareMathOperator{\op}{op}

\usepackage{geometry}
\usepackage{mathtools}

\usepackage{amsmath,amsthm,amssymb,amsfonts,bbm,dsfont}
\usepackage{float}

\usepackage{comment}
\usepackage{bm}
\usepackage{color,soul}
\usepackage{breakcites}     
\usepackage{bigints}

\usepackage{enumitem}

\usepackage{graphicx}
\usepackage{xcolor}
\usepackage{svg}
\usepackage{algorithm}
\usepackage[noend]{algpseudocode}
\usepackage{caption}

\makeatletter

\newcommand{\eps}{\varepsilon }

\newcommand{\F}{\mathrm{F}}

\makeatother

\begin{document}
\title[Neural Entropic Optimal Transport and Gromov-Wasserstein Alignment]{Neural Entropic Optimal Transport and Gromov-Wasserstein Alignment}

\date{Second version: \today}

\author[T. Wang]{Tao Wang}

\address[T. Wang]{Wharton Department of Statistics and Data Science, University of Pennsylvania.}
\email{tawan@wharton.upenn.edu}

\author[Z. Goldfeld]{Ziv Goldfeld}
\address[Z. Goldfeld]{School of Electrical and Computer Engineering, Cornell University.}
\email{goldfeld@cornell.edu}

\begin{abstract}
Optimal transport (OT) and Gromov–Wasserstein (GW) alignment are powerful frameworks for geometrically driven matching of probability distributions, yet their large-scale usage is hampered by high statistical and computational costs. Entropic regularization has emerged as a promising solution, allowing parametric convergence rates via the plug-in estimator, which can be computed using the Sinkhorn algorithm (or its iterations in the GW case). However, Sinkhorn's $O(n^2)$ time complexity for an $n$-sized dataset becomes prohibitive for modern, massive datasets. In this work, we propose a new computational framework for the entropic OT and GW problems that replaces the Sinkhorn step with a neural network trained via backpropagation on mini-batches. By shifting the computational load from the entire dataset to the mini-batch, our approach enables reliable estimation of both the optimal transport/alignment cost and plan at dataset sizes and dimensions far exceeding those tractable with standard Sinkhorn methods. We derive non‐asymptotic error bounds for these estimates, showing they achieve minimax‐optimal parametric convergence rates for compactly supported distributions. Numerical experiments confirm the accuracy of our method in high-dimensional, large-sample regimes where Sinkhorn is infeasible.
\end{abstract}

\thanks{Z. Goldfeld is partially supported by NSF grants CCF-2046018, DMS-2210368, and CCF-2308446, and the IBM Academic Award. Parts of this work were presented at the International Symposium on Information Theory
(ISIT) in 2024 \cite{Wang2024neural}.
}

\maketitle
\section{Introduction}\label{sec:introduction}

\begin{figure}[!t]
    \centering
    \hspace{-4mm}\subfloat[ \textbf{Optimal transport:}\\ \protect{$\mathsf{OT}_c(\mu,\nu)\coloneqq\inf_{\pi\in\Pi(\mu,\nu)}\mathbb{E}_\pi[c(X,Y)]$}, where $c$ is the cost function and the shown $\pi_\text{OT}$ minimizes, e.g., the $p$-Wasserstein cost on $\RR^d$, $c(x,y)=\|x-y\|^p$.]{
        \scalebox{0.5}{\includegraphics{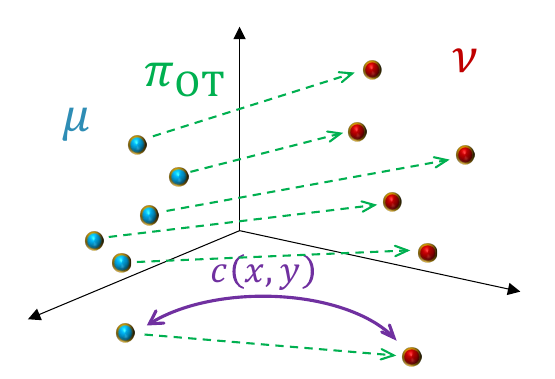}\qquad\qquad\qquad\ \ }
    }
    \qquad
    \subfloat[ \textbf{Gromov-Wasserstein alignment:}\\
    \protect{$\mathsf{GW}_{p, q}(\mu, \nu)^p\coloneqq\inf _{\pi \in \Pi\left(\mu, \nu\right)} \mathbb{E}_{\pi\otimes\pi}\big[\Delta_q\big((X,X'),(Y,Y')\big)^p\big]$}, where $\Delta_q\big((x,x'),(y,y')\big)=\big| \mathsf{d}_\cX(x,x')^q-\mathsf{d}_\cY(y,y')^q\big|$ is the distance distortion cost between the mm spaces $(\cX,\mathsf{d}_\cX,\mu)$ and $(\cY,\mathsf{d}_\cY,\nu)$.]{
        \scalebox{0.5}{\includegraphics{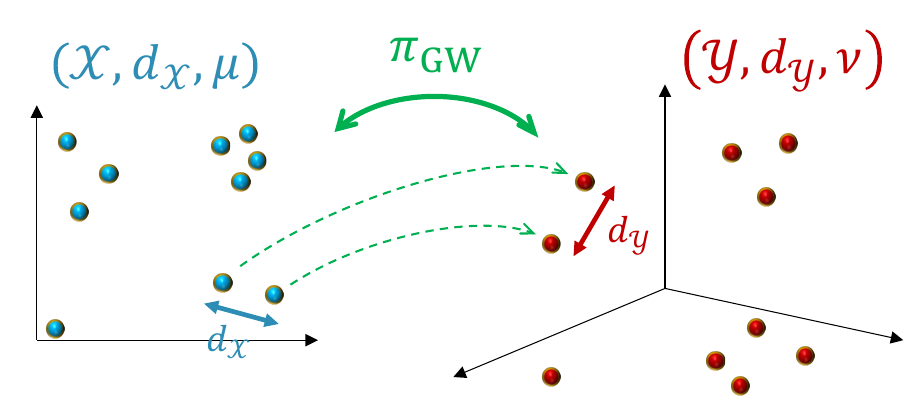}}
    }
    \caption{Illustration of optimal plan for the OT and GW problems between $\mu$ and $\nu$, with $\Pi(\mu,\nu)$ designating the set of all their couplings.}
    \label{fig:OTGW_plans}
\end{figure}

Matching probability distributions in a geometrically meaningful way is a fundamental problem across statistics, probability, and machine learning. Two canonical frameworks to do so are optimal transport (OT), which seeks a coupling that minimizes a prescribed cost function between samples, and Gromov–Wasserstein (GW) alignment, which matches metric measure spaces in a manner that preserves their intrinsic distances; see \cref{fig:OTGW_plans}. Both have proven highly effective in domains encompassing single-cell genomics \cite{schiebinger2017reconstruction,demetci2020gromov,blumberg2020mrec,bunne2024optimal}, computer vision \cite{memoli2009spectral,solomon2015convolutional,li2020continuous,koehl2023computing}, graph matching \cite{chen2020graph,petric2019got,xu2019gromov,xu2019scalable,rioux2024limit}, machine learning \cite{arjovsky2017wasserstein,tolstikhin2017wasserstein,courty2017joint,alvarez2018gromov,bunne2019learning,zhang2024gradient}, and more. However, scaling up these methods to dimensionality and dataset sizes featuring in modern data-driven pipelines is hampered by the high statistical and computational cost of the OT and GW problems. Indeed, empirical convergence rates suffer from the curse of dimensionality \cite{fournier2015rate,zhang2024gromov}, whereby the number of samples needed for reliable estimation grows exponentially with dimension. Computationally, the OT is a linear program solvable in $O(n^3 \log (n))$ time \cite{peyre2017computational}, which is prohibitive for large $n$, while GW alignment is an instance of the quadratic assignment problem, which is NP-hard \cite{commander2005survey}.  

To circumvent the statistical and computational hardness issues, entropic regularization has emerged as a popular remedy, where the original cost is penalized by a  Kullback-Leibler (KL) divergence term, e.g., entropic OT (EOT) reads \cite{cuturi2013sinkhorn}:
\[
\OT_c^{\eps}(\mu,\nu)\coloneqq \inf_{\pi\in\Pi(\mu,\nu)} \EE_\pi[c]+\eps \KL(\pi\|\mu\otimes \nu),
\]
The strongly convex regularizer endows EOT with a unique optimal $\pi_\star^\eps$ solution, leading to $n^{-1/2}$ parametric estimation rates in arbitrary dimension \cite{genevay2019sample,mena2019statistical} and efficient computation via Sinkhorn's algorithm \cite{cuturi2013sinkhorn} in $O(n^2)$ time. For entropic GW (EGW) \cite{solomon2016entropic,peyre2016gromov}, defined analogously to the above but with the regularized distance distortion cost, similar parametric estimation rates apply for Euclidean mm spaces \cite{zhang2024gromov,groppe2023lower}. Recently, \cite{rioux2024entropic} provided the first algorithms with non-asymptotic convergence guarantees for computing EGW via accelerated first-order methods with a Sinkhorn oracle. In both EOT and EGW cases, however, the $O(n^2)$ runtime of Sinkhorn's algorithm quickly becomes intractable for massive datasets that appear in modern machine learning pipelines. Scaling up these optimal matching frameworks to regimes where $n \gg 1$ is prohibitively large is the main motivation of this work.

\subsection{Contributions} \label{subsec:Contributions}

This work develops a provably accurate neural estimation framework for EOT and EGW. Our estimators are trained end-to-end via backpropagation on mini-batches, effectively shifting the computational burden from the full dataset size $n$ to the typically much smaller mini-batch size $b$. We begin with EOT neural estimation, which forms the basis for our subsequent EGW estimator. Let $c$ be a general (smooth) cost function, and recall the semi-dual form
\[
\OT_c^\eps(\mu,\nu)=
  \sup_{\varphi\in L^1(\mu)} 
      \EE_\mu[\varphi] + \EE_\nu[\varphi^{c,\eps}],
\]
where $\varphi^{c,\eps}$ is the $(c,\eps)$-transform of $\varphi$ with respect to (w.r.t.) the cost function $c$ (see \eqref{eq:EOT_semidual}). Our \emph{EOT neural estimator} (NE) is obtained by parameterizing the dual potential with a neural network (NN), approximating expectations via sample means, and optimizing the resulting empirical objective over the NN parameters. The approach yields not only an estimate of the EOT cost, but also of the corresponding coupling, as the EOT coupling and optimal potential are related through the Gibbs density $\frac{d\pi^\eps_{\star}}{d \mu \otimes \nu}(x,y)=e^{\frac{\varphi(x)+\varphi^{c,\eps}(y)-c(x,y)}{\eps}}$. The NE thus serves as an alternative to the standard plug-in-plus-Sinkhorn approach, offering a scalable, data-driven proxy for both the EOT cost and plan.

Building on the EOT NE, we move to the more challenging task of estimating the EGW distance. We focus on the nominal case of the quadratic EGW distance, i.e., when $p=q=2$ (these order indices are henceforth omitted from our notation), between Euclidean mm spaces $\left(\mathbb{R}^{d_x},\|\cdot\|, \mu\right)$ and $\left(\mathbb{R}^{d_y},\|\cdot\|, \nu\right)$. Thanks to the recently developed EGW variational representation from \cite[Theorem 1]{zhang2024gromov}, we have
\begin{equation}
\mathsf{GW}^{\eps}(\mu, \nu)
   = \mathsf{C}_{\mu, \nu}+\inf _{\bA  \in \mathcal{D}_{M_{\mu,\nu}}}\big\{32\|\bA \|_\F^2+\mathsf{OT}^\eps_\bA(\mu,\nu)\big\},\label{EQ:dual_intro}
\end{equation}
where $\mathsf{C}_{\mu, \nu}$ is a constant that depends only on the marginals (and is easy to estimate; see \eqref{eq:C_munu_constant} ahead),  $\mathcal{D}_{M_{\mu,\nu}} \subsetneq \mathbb{R}^{d_x \times d_y}$ is compact whenever $\mu,\nu$ have finite 2nd moments, while $\mathsf{OT}^\eps_\bA$ is the EOT problem with cost function $c_\bA(x,y)= -4\|x\|^2\|y\|^2-32x^{\intercal}\bA  y$. As our neural EOT estimator can treat any smooth cost $c$, we propose to estimate EGW cost and plan by replacing $\mathsf{OT}^\eps_\bA$ with its NE and jointly optimize over the NN parameters and $\bA$ via a gradient-based routine. The computational cost of the EGW NE again scales with the mini-batch size $b$ (as opposed to $n$), resulting in a marked increase in scalability compared to Sinkhorn. In addition, the neural EOT and EGW estimators can be seamlessly integrated into larger machine learning pipelines, e.g., as a loss or regularizer.

For both EOT and EGW, we provide formal guarantees on the quality of the NE of the cost and optimal plan. Our analysis relies on non-asymptotic function approximation theorems and tools from empirical process theory to bound the two sources of error involved: function approximation and empirical estimation. Given $n$ samples from the populations, we show that in both cases the effective error of a NE realized by a shallow NN of $k$ neurons scales as
    \begin{equation}
    \label{effictive_error}
    O_{d_x,d_y}\left(\operatorname{poly}(1 / \varepsilon)\left(k^{-1 / 2}+n^{-1 / 2}\right)\right),
    \end{equation}
    where the subscript on the $O(\cdot)$ notation indicates that the hidden constants depend on the data dimensions $d_x$ and $d_y$. Crucially, the degree of the polynomial in $1 / \varepsilon$ also scales with the dimension, with the precise form of this dependence is explicitly characterized in our main results (see Theorems \ref{thm:Bound_EOT1}-\ref{thm:Bound_NeuralCoupling_EOT} and \ref{thm:Bound_EGW1}-\ref{thm:Bound_NeuralCoupling}).
The bounds on the EOT and EGW cost estimation errors hold for arbitrary, compactly supported distributions. This stands in stark contrast to existing neural estimation error bounds for other divergences \cite{nguyen2010estimating,alvarez2018gromov,sreekumar2022neural,goldfeld2022k,tsur2023max}, which typically require strong regularity assumptions on the population distributions (e.g., H\"older smoothness of densities). This is unnecessary in our setting thanks to the inherent regularity of dual EOT potentials for smooth cost functions.

The above bound reveals the optimal scaling of the NN and dataset sizes, namely $k \asymp n$, which achieves the parametric convergence rate of $n^{-1/2}$ and guarantees minimax-rate optimality of our NE for compactly supported distributions. The explicitly characterized polynomial dependence on $1/\eps$ in our bound is in line with the bounds for EOT and EGW estimation via empirical plug-in \cite{mena2019statistical,zhang2024gromov,groppe2023lower}. We also note that our neural estimation results readily extend to the inner product EGW distance, which has seen notable interest due to its analytic tractability \cite{dumont2022existence,le2022entropic,zhang2024gradient}. We empirically test the developed NE on synthetic and real-world datasets, demonstrating its scalability to high dimensions and validating our theory.

\subsection{Related Literature}\label{subsec:lietrature_review}

Neural estimation is a popular approach for enhancing scalability. Prior research explored the tradeoffs between approximation and estimation errors in non-parametric regression \cite{barron1994approximation,bach2017breaking,suzuki2018adaptivity} and density estimation \cite{yang1999information,uppal2019nonparametric} tasks. More recently, neural estimation of statistical divergences and information measures has been gaining attention. The mutual information NE (MINE) was proposed in \cite{belghazi2018mine}, and has seen various improvements since \cite{poole2018variational,song2019understanding,chan2019neural,mroueh2021improved}. Extensions of the neural estimation approach to directed information were studied in \cite{molavipour2021neural,tsur2023neural,tsur2023data}. Theoretical guarantees for $f$-divergence NEs, accounting for approximation and estimation errors, as we do here, were developed in \cite{sreekumar2021non,sreekumar2022neural} (see also \cite{nguyen2010estimating} for a related approach based on reproducing kernel Hilbert space parameterization). Neural estimation of the Stein discrepancy and the minimum mean squared error were considered in \cite{repasky2023neural} and \cite{diaz2021lower}, respectively. Neural methods for approximate computation  of the Wasserstein distances have been considered under the Wasserstein generative adversarial network (GAN) framework \cite{arjovsky2017wasserstein,gulrajani2017improved}, although these approaches lack formal guarantees. A neural computation framework for entropically and quadratically regularized OT was proposed in \cite{seguy2018large}, also without guarantees and while optimizing over two dual potentials (as opposed to our approach which computes a single potential, thanks to the semi-dual EOT form). More recently, \cite{daniels2021score} studied a score-based generative neural EOT model, while an energy-based model was considered in \cite{mokrov2023energy}.

Computationally tractable relaxation and reformulations of the GW problem has seen significant interest. The sliced GW distance \cite{vayer2019sliced} seeks to reduce the computational burden by averaging GW costs between one-dimensional projections of the marginal distributions. The utility of this approach, however, is contingent on resolving the GW problem on $\RR$, which remains open \cite{beinert2023assignment}. The unbalanced GW distance was proposed in \cite{sejourne2021unbalanced}, along with a computationally tractable convex relaxation thereof. A Gromov-Monge map-type formulation was explored in \cite{zhang2022cycle}, which directly optimizes over bi-directional maps to attain structured solutions. Yet, it is the entropically regularized GW distance \cite{peyre2016gromov,solomon2016entropic} that has been widely adopted in practice, thanks to its compatibility with iterative methods based on Sinkhorn's algorithm \cite{cuturi2013sinkhorn}. Under certain low-rank assumptions on the cost matrix, \cite{scetbon2022linear} presented an adaptation of the mirror descent approach from \cite{peyre2016gromov} that speeds up its runtime from $O(n^3)$ to $O(n^2)$. More recently, \cite{rioux2024entropic} proposed an accelerated first-order method based on the dual formulation from \eqref{EQ:dual_intro} and derived non-asymptotic convergence guarantees for it. To the best of our knowledge, the latter is the only algorithm for computing a GW variant subject to formal convergence claims. However, all the above EGW methods use iterates that run Sinkhorn's algorithm, whose time and memory complexity scale quadratically with $n$. This hinders applicability to large-scale problems.

\section{Background and Preliminaries}\label{sec:background}
\subsection{Notation}\label{subsec:Notation}
Let $\|\cdot\|$ denote the Euclidean norm and $\langle \cdot, \cdot \rangle$ designate the inner product. The Euclidean ball is designated as $B_d(r)\coloneqq \{x\in \RR^d:\|x\|<r\}$. We use $\|\cdot\|_{\op}$ and $\|\cdot\|_\F$ for the operator and Frobenius norms of matrices, respectively. For $1 \leq p < \infty$, the $L^p$ space over $\mathcal{X} \subseteq \mathbb{R}^d$ with respect to (w.r.t.) the measure $\mu$ is denoted by $L^p(\mu)$, with $\|\cdot\|_{p, \mu}$ representing the norm. For $p=\infty$, we use $\|\cdot\|_{\infty, \mathcal{X}}$ for standard sup-norm on $\mathcal{X} \subseteq \mathbb{R}^d$. Slightly abusing notation, we also set $\|\mathcal{X}\|\coloneqq\sup _{x \in \mathcal{X}}\|x\|_{\infty}$.

The class of Borel probability measures on $\mathcal{X} \subseteq \mathbb{R}^d$ is denoted by $\mathcal{P}(\mathcal{X})$. To stress that the expectation of $f$ is taken w.r.t. $\mu \in \mathcal{P}(\mathcal{X})$, we write $\mathbb{E}_\mu[f]\coloneqq\int f d \mu$. For $\mu, \nu \in \mathcal{P}(\mathcal{X})$ with $\mu \ll \nu$, i.e., $\mu$ is absolutely continuous w.r.t. $\nu$, we use $\frac{d\mu}{d\nu}$ for the Radon-Nikodym derivative of $\mu$ w.r.t. $\nu$. The subset of probability measures that are absolutely continuous w.r.t. Lebesgue is denoted by $\mathcal{P}_{\mathsf{ac}}(\cX)$. For $p\in[1, \infty)$, further let $\cP_p(\cX)\subseteq \mathcal{P}(\mathcal{X})$ contain only measures with finite $p$-th absolute moment, i.e., $M_p(\mu)\coloneqq \int_{\cX}\|x\|^p d \mu(x)<\infty$ for any $\mu\in\cP_p(\cX)$. 

For any multi-index $\alpha=(\alpha_1,\dots,\alpha_d) \in \mathbb{N}_0^d$ with $|\alpha| = \sum_{j=1}^ d \alpha_j$ ($\NN_0 = \NN \cup \{ 0 \}$), define the differential operator $D^\alpha = \frac{\partial^{|\alpha|}}{\partial x_1^{\alpha_1} \cdots \partial x_{d}^{\alpha{}_d}}$
with $D^0 f = f$. We write $N(\delta, \cF, \mathsf{d})$ for the $\delta$-covering number of a function class $\cF$ w.r.t. a metric $\mathsf{d}$, and $N_{[\,]}(\delta, \cF, \mathsf{d})$ for the bracketing number. For an open set $\mathcal{U} \subseteq \mathbb{R}^d$, $b \geq 0$, and an integer $m \geq 0$, let $\cC_b^m(\mathcal{U})\coloneqq \big\{f \in \cC^m(\mathcal{U}):\,\max _{\alpha:|\alpha| \leq m}\left\|D^\alpha f\right\|_{\infty, \mathcal{U}} \leq b\big\}$ denote the H\"older space of smoothness index $m$ and radius $b$. The restriction of $f: \mathbb{R}^d \rightarrow \mathbb{R}$ to a subset $\mathcal{X} \subseteq \mathbb{R}^d$ is denoted by $f\big|_{\mathcal{X}}$. We use $\lesssim_x$ to denote inequalities up to constants that only depend on $x$; the subscript is dropped when the constant is universal. For $a,b \in \R$, we use the shorthands $a \vee b = \max \{ a,b \}$ and $a \land b = \min \{ a,b \}$.

\subsection{Entropic Optimal Transport}\label{subsec:EOT}

Given distributions $(\mu,\nu)\in\cP(\cX)\times \cP(\cY)$ and a cost function $c:\cX\times \cY\to\RR$, the primal EOT formulation is a convexification of classical OT by means of a KL divergence penalty:
\begin{equation}
\label{EOT}
      \OT_c^{\eps}(\mu,\nu)\coloneqq \inf_{\gamma\in\Pi(\mu,\nu)} \mathbb{E}_\gamma[c]+\eps \KL(\gamma\|\mu\otimes \nu),
\end{equation}
where $\eps>0$ is a regularization parameter and $\KL(\mu\|\nu)\coloneqq\mathbb{E}_\mu\left[\log \left(\frac{d \mu}{d \nu}\right) \right]$ if $\mu \ll \nu$ and $+\infty$ otherwise. Classical OT \cite{villani2008optimal,santambrogio2010} is obtained from \eqref{EOT} by setting $\eps=0$.\footnote{We use $\gamma\in\Pi(\mu,\nu)$ for the EOT plan to differentiate it from the GW alignment plan, for which we reserve the symbol $\pi$.} When~$c\in L^1(\mu\otimes \nu)$, EOT admits the dual and semi-dual formulations, which are, respectively, given~by
\begin{align}
\OT_c^{\eps}(\mu,\nu)&=\sup_{(\varphi,\psi)\in L^1(\mu)\times L^1(\nu)}\mathbb{E}_\mu[\varphi] +\mathbb{E}_\nu [\psi]-\eps \mathbb{E}_{\mu\otimes\nu}\left[e^{\frac{\varphi\oplus \psi-c}{\eps}}\right] +\eps, \label{eq:EOT_dual}\\
&=\sup_{\varphi\in L^1(\mu)}\mathbb{E}_\mu[\varphi]+\mathbb{E}_\nu[\varphi^{c,\eps}],\label{eq:EOT_semidual}
\end{align}
where we have defined $(\varphi\oplus \psi)(x,y)=\varphi(x)+\psi(y)$ and the $(c,\eps)$-transform of $\varphi$ is given by $\varphi^{c,\eps}=-\eps \log \left(\int_{\mathcal{X}} \exp \left(\frac{\varphi(x)-c(x,\cdot)}{\eps}\right)d\mu\right)$. There exist functions $(\varphi,\psi)$ that achieve the supremum in $\eqref{eq:EOT_dual}$, which we call \emph{EOT potentials}. These potentials are almost surely (a.s.) unique up to additive constants, i.e., if $(\tilde{\varphi}, \tilde{\psi})$ is another pair of EOT potentials, then there exists a constant $a \in \mathbb{R}$ such that $\tilde{\varphi}=\varphi+a$ $\mu$-a.s. and $\tilde{\psi}=\psi-a$ $\nu$-a.s. 

\medskip
A pair $(\varphi, \psi) \in L^1(\mu) \times L^1(\nu)$ are EOT potentials if and only if they satisfy the Schr\"odinger system
\begin{equation}
    \label{eq:Schrodinger system}
\int e^{\frac{\varphi(x)+\psi(\cdot)-c(x, \cdot)}{\eps}} d \mu(x)=1 \quad \nu \text {-a.s. } \quad \text { and } \quad \int e^{\frac{\varphi(\cdot)+\psi(y)-c(\cdot, y)}{\eps}} d \nu(y)=1 \quad \mu \text {-a.s. }
\end{equation}
Furthermore, $\varphi$ solves the semi-dual from \eqref{eq:EOT_semidual} if an only if $(\varphi,\varphi^{c,\eps})$ is a solution to the full dual in \eqref{eq:EOT_dual}. Given EOT potentials $(\varphi, \psi)$, the unique EOT plan can be expressed in their terms as
\begin{equation}
d\gamma^\eps_{\star}=e^{\frac{\varphi\oplus \psi-c}{\eps}} d \mu \otimes \nu.\label{eq:EOT_plan}
\end{equation}
Subject to smoothness assumptions on the cost function and the population distributions, various regularity properties of EOT potentials can be derived; cf., e.g., \cite[Lemma 1]{goldfeld2022limit}.    

\subsection{Entropic Gromov-Wasserstein Distance}\label{subsec:EGW} 

We consider the quadratic EGW distance (i.e., when $p=q=2$) between Euclidean mm spaces $(\RR^{d_x},\|\cdot\|,\mu)$ and $(\RR^{d_y},\|\cdot\|,\nu)$, where $\mu,\nu$ are assumed to have finite 4th absolute moments:
\begin{equation}
\label{eq:QuadEGW}
    \mathsf{GW}^{\eps}(\mu,\nu)=\inf_{\pi\in\Pi(\mu,\nu)}\mathbb{E}_{\pi\otimes\pi}\left[ \left|\|X-X'\|^2-\|Y-Y'\|^2\right|^2\right]+\eps \KL(\pi\|\mu\otimes \nu). 
\end{equation}
By a standard compactness argument, the infimum above is always achieved, and we call such a solution an \emph{EGW alignment plan}, denoted by $\pi_\star^\eps$ \cite{memoli2011gromov}. By analogy to OT, EGW serves as a proxy of the standard GW distance up to an additive gap of $O\big(\eps\log(1/\eps)\big)$ \cite{zhang2024gromov}. One readily verifies that, like the unregularized distance,~EGW is invariant to isometric actions on the marginal spaces such as orthogonal rotations and translations.

\medskip
It was shown in \cite[Theorem 1]{zhang2024gromov} that when $\mu,\nu$ are centered, which we may assume without loss of generality, the EGW cost admits a dual representation. To state it, define 
\begin{equation}
    \mathsf{C}_{\mu,\nu}\coloneqq\mathbb{E}_{\mu\otimes\mu}\big[\|X-X'\|^4\big]+\mathbb{E}_{\nu\otimes\nu}\big[\|Y-Y'\|^4\big]-4\mathbb{E}_{\mu\otimes\nu}\big[\|X\|^2\|Y\|^2\big],\label{eq:C_munu_constant}
\end{equation}
which, evidently, depends only on the marginals $\mu,\nu$. We have the following.

\begin{lemma}[EGW duality; Theorem 1 in \cite{zhang2024gromov}]
\label{lemma:EGWDual}
    Fix $\eps>0$, $(\mu,\nu)\in\cP_4(\RR^{d_x})\times\cP_4(\RR^{d_y})$ with zero mean, and let $M_{\mu,\nu}\coloneqq \sqrt{M_2(\mu)M_2(\nu)}$. Then, 
    \begin{equation}
    \label{eq:S2Decomp}
       \mathsf{GW}^{\eps}(\mu,\nu)=\sC_{\mu,\nu}+\inf_{\bA \in\RR^{d_x\times d_y}}\Big\{32\|\bA \|_\F^2+\OT_\bA ^{\mspace{1mu}\eps}(\mu,\nu)\Big\},
     \end{equation}
    where $\OT_\bA ^{\mspace{1mu}\eps}$ is the EOT problem with cost $c_{\bA}:(x,y)\in\RR^{d_x}\times \RR^{d_y}\mapsto-4\|x\|^2\|y\|^2-32x^{\intercal}\bA  y$. 
    Moreover, the infimum is achieved at some $\bA _{\star}\in\cD_{M_{\mu,\nu}}\coloneqq [-M_{\mu,\nu}/2,M_{\mu,\nu}/2]^{d_x\times d_y}$. 
\end{lemma}

Although \eqref{eq:S2Decomp} illustrates a connection between the EGW and EOT problems, the outer minimization over $\bA $ necessitates studying EOT with an \textit{a priori} unknown cost function $c_{\bA}$. To enable that, in \cref{lemma:regularity_EOT} ahead we show that there exist smooth EOT potentials $(\varphi_{\bA},\varphi_{\bA}^{c,\eps})$ for $\OT_\bA ^{\mspace{1mu}\eps}(\mu,\nu)$ satisfying certain derivative estimates uniformly in $\bA \in\mathcal D_M$ (for some $M\geq M_{\mu,\nu}$) and $\mu,\nu$. This enables accurately approximating them by NNs and, in turn, unlocks our neural estimation approach.

\begin{remark}[Inner product cost]
\label{rem:EIGW}
Our NE approach for the quadratic EGW distance extends almost directly to the case of the inner product distortion cost (abbreviated henceforth by EIGW)
    \begin{equation}  \mathsf{IGW}^{\eps}(\mu,\nu)=\inf_{\pi\in\Pi(\mu,\nu)}\mathbb{E}_{\pi\otimes\pi}\left[ \left|\langle X,X'\rangle-\langle Y,Y'\rangle\right|^2\right]+\eps \KL(\pi\|\mu\otimes \nu).\label{eq:EIGW}   
\end{equation}
The cost function above is not a distance distortion measure; rather, it quantifies the change in angles. This object has received recent attention due to its analytic tractability \cite{le2022entropic} and since it captures a meaningful notion of discrepancy between mm spaces with a natural inner product structure. EIGW enjoys a similar decomposition, dual representation, and regularity of dual potentials as its quadratic counterpart. As such, it falls under our neural estimation framework and can be treated similarly. To avoid repetition, we collect the results pertaining to the EIGW distance in \cref{appen:EIGW}. 
\end{remark}

\section{Neural Estimation of EOT Cost and Plan}


 We provide non-asymptotic error bounds for the NE of EOT, accounting both for the cost and optimal plan. Generally, the proposed approach entails three sources of error: (i) function approximation of the dual EOT potentials by neural nets; (ii) empirical estimation of the means by sample averages; and (iii) optimization, which comes from the suboptimality of gradient-based routines employed in practice. We provide sharp bounds on the errors of types (i) and (ii). Treating error type (iii) would require global optimality guarantees for the algorithm employed to solve NE of EOT, which is an instance of the general open question of convergence analysis for nonlinear neural network optimization over nonconvex loss landscapes. Still, we provide a partial account of the optimization error in two ways: (i) demonstrate that it decomposes out of the error analysis (see \cref{rem:opt_error}), resulting in an another additive term to the bounds presented in Theorems \ref{thm:Bound_EOT1}-\ref{thm:Bound_EOT_2} and \ref{thm:Bound_EGW1}-\ref{thm:Bound_EGW2} ahead; and (ii) establish the convergence of our gradient-based algorithm (\cref{alg:EGW_NE}) given the inner neural net optimization can be solved into  global optimum up to error $\delta^\prime$, see \cref{subsec:experiments_algorithm} and \cref{lem:delta-oracle}.

We consider compactly supported distributions $(\mu,\nu)\in\cP(\cX)\times\cP(\cY)$, and assume, for simplicity, that $\cX\subseteq [-1,1]^{d_x}$ and $\cY\subseteq [-1,1]^{d_y}$ (our results readily extend to arbitrary compact supports in Euclidean mm spaces). Without loss of generality, further suppose that $d_x\leq d_y$, and consider a general smooth cost function $c: \RR^{d_x}\times \RR^{d_y} \rightarrow \mathbb{R}$. We next describe the NE for the EOT distance and provide non-asymptotic performance guarantees for estimation of both the transport cost and plan. All proofs are deferred to the \cref{appen:Proofs}.

\subsection{EOT Neural Estimator}
\label{subsec:EOT_Neural_Estimator}
Let $X^n\coloneqq\left(X_1, \cdots, X_n\right)$ and $Y^n\coloneqq\left(Y_1, \cdots, Y_n\right)$ be $n$ independently and identically distributed (i.i.d.) samples from $\mu$ and $\nu$, respectively. Further suppose that the sample sets are independent of each other. Denote the corresponding empirical measures by $\hat\mu_n=n^{-1}\sum_{i=1}^n\delta_{X_i}$ and $\hat\nu_n=n^{-1}\sum_{i=1}^n\delta_{Y_i}$. Our NE is realized by a shallow ReLU NN (i.e., with a single hidden layer) with $k$ neurons, which defines the function class
\begin{equation}
\label{NN_class}
    \cF_{k,d}(a)\coloneqq \left\{f: \mathbb{R}^{d} \rightarrow \mathbb{R}:\begin{aligned}
&f(x)=\sum_{i=1}^k \beta_i \phi\left(w_i \cdot x+b_i\right)+w_0 \cdot x+b_0, \\
& \max _{1 \leq i \leq k}\left\|w_i\right\|_1 \vee\left|b_i\right| \leq 1,\ \max _{1 \leq i \leq k}\left|\beta_i\right| \leq 2ak^{-1},\ \left|b_0\right| \leq a,\ \left\|w_0\right\|_1 \leq a
\end{aligned}\right\},
\end{equation}
where $a\in\RR_{\geq 0}$ specifies the parameter bounds and $\phi:\RR\to\RR_{\geq 0}: z\mapsto z\vee 0$ is the ReLU activation function, which acts on vectors component-wise.

We parametrize the semi-dual form of $\OT_c^{\mspace{1mu}\eps}(\mu,\nu)$ (see \eqref{eq:EOT_semidual}) using a NN from the class $\cF_{k,d}(a)$ and replace expectations with sample means. Specifically, the EOT distance NE is
\begin{equation}\label{eq:NE_EOTc}
    \widehat\OT_{c, k, a}^{\mspace{1mu}\eps}(X^n, Y^n)\coloneqq 
   \max_{f \in \cF_{k, d_x}(a)} \frac{1}{n}\sum_{i=1}^nf(X_i)-\frac{\eps}{n}\sum_{j=1}^n\log\left(\frac{1}{n}\sum_{i=1}^n\exp\left(\frac{f(X_i)-c(X_i,Y_j)}{\eps}\right)\right). 
\end{equation}
The objective above can be optimized using mini-batch gradient descent over NN the parameters. Unlike the Sinkhorn algorithm, which must process the full $n$-sized dataset and suffers from scalability issues when $n$ is prohibitively large, our neural estimation approach only accesses the dataset through mini-batches. The NE objective is calculated on a mini-batch of size $b$ in $O(b)$ time, which amounts to $O(n)$ for a single epoch. This enables running the NE on large-scale problems, beyond the reach of Sinkhorn's algorithm, as illustrated in \cref{sec:experiments}. Upon convergence, the optimized neural dual potential $f_\star\in\cF_{k,d_x}(a)$ gives rise to a neural EOT plan
\begin{equation}
    d \gamma^\eps_{f_\star}(x, y) \coloneqq \frac{\exp \left(\frac{f_\star(x)-c(x,y)}{\eps}\right)}{\int_{\cX} \exp \left(\frac{f_\star(x)-c(x,y)}{\eps}\right) d\mu(x)} d\mu\otimes\nu(x,y).\label{eq:NE_EOTc_plan}
\end{equation}
which serves as an estimate of the true EOT plan $\gamma^\eps_\star$ from \eqref{eq:EOT_plan}.

\subsection{Performance Guarantees}
We provide formal guarantees for the neural estimator of the EOT cost and plan defined above. Starting from the cost estimation setting, we establish two separate bounds on the effective (approximation plus estimation) error. The first is non-asymptotic and presents optimal convergence rates, but calibrates the NN parameters to a cumbersome dimension-dependent constant. Following that, we present an alternative bound that avoids the dependence on the implicit constant, but at the expense of a polylogarithmic slow-down in the rate and a requirement that the NN size $k$ is large enough.

\begin{theorem}[EOT cost neural estimation; bound 1] \label{thm:Bound_EOT1}
     There exists a constant $C>0$ depending only on $c, d_x, d_y$, such that setting $a=C(1+\eps^{1-s})$ with $s=\left\lfloor d_x  / 2\right\rfloor+3$, we have
\begin{equation}
\label{Bound_EOT1}
\begin{aligned}
&\sup _{(\mu, \nu) \in \mathcal{P}(\cX) \times \mathcal{P}(\cY)} \mathbb{E}\left[\left|\widehat\OT_{c,k,a}^{\mspace{1mu}\eps}\left(X^n, Y^n\right)-\OT_c^{\eps}(\mu, \nu)\right|\right]\\
&\lesssim_{c, d_x, d_y}\left(1+\frac{1}{\eps^{\left\lfloor \frac{d_x}{2}\right\rfloor+2}}\right) k^{-\frac{1}{2}}+\min\left\{1+\frac{1}{\eps^{\left\lceil d_x+\frac{d_y}{2}\right\rceil+4}}\,,\left(1+\frac{1}{\eps^{\left\lfloor \frac{d_x}{2}\right\rfloor+2}}\right)\sqrt{k}\right\} n^{-\frac{1}{2}},
\end{aligned}
\end{equation}
where the dependence on the cost function $c$ is through $\|c\|_\infty$ and bounds on its derivatives. 
\end{theorem}

\begin{proof}
To analyze the neural estimator of the EOT cost, we decompose the error into the approximation and empirical estimation errors:
\begin{align}
&\mathbb{E}\left[\left|\widehat\OT_{c,k,a}^{\mspace{1mu}\eps}\left(X^n, Y^n\right)-\OT_c^{\eps}(\mu, \nu)\right|\right]\nonumber\\
&\qquad\qquad\leq\underbrace{\left|\OT_{c,k,a}^{\mspace{1mu}\eps}(\mu, \nu)-\OT_c^{\eps}(\mu, \nu)\right|}_{\text{ Approximation error }}+\underbrace{\mathbb{E}\left[\left|\OT^{\eps}_{c,k,a}(\mu, \nu)-\widehat{\OT}^{\eps}_{c,k,a}\left(X^n, Y^n\right)\right|\right]}_{\text{ Estimation error }},\label{eq:EOT_NE_error_decomposition}
\end{align}
where the population-level neural EOT cost $\OT^{\eps}_{c,k,a}(\mu, \nu)$ is defined as
\begin{equation}
\OT^{\eps}_{c,k,a}(\mu, \nu)\coloneqq\sup_{f\in \cF_{k,d_x}(a)}\int f d\mu+\int f^{c,\eps}d\nu.\label{eq:neural_pop_EOT}
\end{equation}
We analyze each term separately and summarize the results in the following technical lemmas.

\begin{lemma}[Approximation error bound]
\label{lem:approximation_error_bound}
Under the setting of \cref{thm:Bound_EOT1}, we have

\begin{equation}
\label{Approximation_Bound_EOT}
\left|\OT_{c, k, a}^{\eps}(\mu, \nu)-\OT_c^{\eps}(\mu, \nu)\right| \lesssim_{c, d_x, d_y}\left(1+\frac{1}{\eps^{\left\lfloor\frac{d_x}{2}\right\rfloor}+2}\right) k^{-\frac{1}{2}}
\end{equation}
\end{lemma}
The proof of \cref{lem:approximation_error_bound} is given in \cref{appen:proof_approx_error_bound}. We establish  regularity of semi-dual EOT potentials (namely, $(\varphi,\varphi^{c,\eps})$ in \eqref{eq:EOT_semidual}), showing that they belong to a H\"older class of arbitrary smoothness. This, in turn, allows accurately approximating these dual potentials by NNs from the class $\cF_{k,d}(a)$ with error $O(k^{-1/2})$, yielding the approximation bound.

\begin{lemma}[Estimation error]
\label{lem:estimation_error_bound}
Under the setting of \cref{thm:Bound_EOT1}, we have 
\begin{equation}
\label{Estimation_Bound_EOT}
\begin{aligned}
    \mathbb{E}&\left[\left|\OT^\eps_{c ,k,a}(\mu, \nu)-\widehat\OT^\eps_{c , k,a}\left(X^n, Y^n\right)\right|\right]\\
    &\qquad\qquad\qquad\qquad\qquad\lesssim_{c, d_x, d_y} \min\left\{1+\frac{1}{\eps^{\left\lceil d_x+\frac{d_y}{2}\right\rceil+4}}\,,\left(1+\frac{1}{\eps^{\left\lfloor \frac{d_x}{2}\right\rfloor+2}}\right)\sqrt{k}\right\} n^{-\frac{1}{2}}.
\end{aligned}
\end{equation}
\end{lemma}
The proof of \cref{lem:estimation_error_bound}, given in \cref{appen:proof_estimation_error}, employs standard maximal inequalities from empirical process theory along with a bound on the covering or bracketing number of $(c,\eps)$-transform of the NN class. To improve the dependence of the bound on dimension, we also leverage the lower complexity adaptation (LCA) principle from \cite{hundrieser2022empirical}.  
\medskip

Inserting the bounds from Lemmas \ref{lem:approximation_error_bound} and \ref{lem:estimation_error_bound} into \eqref{eq:EOT_NE_error_decomposition} yields the result of \cref{thm:Bound_EOT1}.
\end{proof}

\begin{remark}[Minimax optimality]
Evidently, by equating $n\asymp k$, the bound from \cref{thm:Bound_EOT1} yields the $n^{-1/2}$ parametric convergence rate, which is sharp. This implies minimax optimality of our NE for the EOT cost over the class of compactly supported distributions as above. 
\end{remark}

\begin{remark}[Almost explicit expression for $C$]\label{REM:C_constant}
The expression of the constant $C$ in \cref{thm:Bound_EOT1} is cumbersome, but can nonetheless be evaluated. Indeed, one may express $C=C_sC_{c, d_x, \mathcal{X}, \mathcal{Y}} \bar{c}_{d_x}$, with explicit expressions for $C_{c, d_x, \mathcal{X}, \mathcal{Y}}$ and $\bar{c}_{d_x}$ given in \eqref{eq:constant_C_dx_dy}
and \eqref{eq:constant_c_bar}, respectively, while $C_s$ is a combinatorial constant that arises from
the multivariate Faa di Bruno formula (cf. \eqref{eq:constant_C_s_1}-\eqref{eq:constant_C_s_2}). The latter constant is quite convoluted and is the main reason we view $C$ as implicit. 
\end{remark}

Our next bound circumvents the dependence on $C$ by letting the NN parameters grow with its size $k$. This bound, however, requires $k$ to be large enough (specifically when k is such that $m_k$ is larger than $C\left(1+\varepsilon^{1-s}\right)$ given in \cref{thm:Bound_EOT1}) and entails additional polylog factors in the rate. The proof is similar to that of \cref{thm:Bound_EOT1} and given in \cref{appen:proof_Bound_EOT2}.
\begin{theorem}[EOT cost neural estimation; bound 2]
    \label{thm:Bound_EOT_2}
       Let $\epsilon>0$ and set $m_k=\log k\vee 1$. Assuming $k$ is sufficiently large, we have
       \begin{equation}
           \label{Bound_EOT2}
               \begin{aligned}
   &\sup _{(\mu, \nu) \in \mathcal{P}(\cX) \times \mathcal{P}(\cY)} \mspace{-5mu}\mathbb{E}\left[\left|\widehat\OT_{c,k,a}^{\mspace{1mu}\eps}\left(X^n, Y^n\right)-\OT_c^{\eps}(\mu, \nu)\right|\right]\\
   &\qquad\qquad\lesssim_{c, d_x, d_y} \left(1+\frac{1}{\eps^{\left\lfloor\frac{d_x}{2}\right\rfloor+2}}\right) k^{-\frac{1}{2}}+\min\left\{\left(1+\frac{1}{\eps^{\left[\frac{d_y}{2}\right\rceil}}\right)(\log k)^2\,, \sqrt{k}\log k\right\} n^{-\frac{1}{2}}.
\end{aligned}
       \end{equation}
\end{theorem}

Lastly, we account for the quality of the neural plan estimate from  \eqref{eq:NE_EOTc_plan} by comparing it, in KL divergence, to the true EOT plan $\gamma_{\star}^{\varepsilon}$.
\begin{theorem}[EOT plan neural estimation]
\label{thm:Bound_NeuralCoupling_EOT}
Suppose that $\mu \in \mathcal{P}_{\mathrm{ac}}(\mathcal{X})$. Let $\hat{f}_{\star}$ be a maximizer of $\widehat{\mathrm{OT}}_{c,k, a}^{\varepsilon}\left(X^n, Y^n\right)$ from \eqref{eq:NE_EOTc}, with a as defined in \cref{thm:Bound_EOT1}. Then, the induced neural plan $\gamma_{\hat{f}_*}^{\varepsilon}$ from \eqref{eq:NE_EOTc_plan} satisfies
\begin{equation}
\mathbb{E}\left[\KL\left(\gamma_{\star}^{\varepsilon} \middle\| \gamma_{\hat{f}_{\star}}^{\varepsilon}\right)\right]\lesssim_{c,d_x,d_y}\varepsilon^{-1} \delta_{\mathrm{EOT}}
\end{equation}
where $\delta_{\mathrm{EOT}}$ denotes the right-hand side (RHS) bound of \eqref{Bound_EOT1} in \cref{thm:Bound_EOT1}.
\end{theorem}

\cref{thm:Bound_NeuralCoupling_EOT} is proved in \ref{appen:proof_NeuralCoupling_EOT}. The key step in the derivation shows that the KL divergence between the alignment plans, in fact, equals the gap between the EOT cost $\OT^{\varepsilon}_c$ and its neural estimate from \eqref{eq:NE_EOTc}, up to a multiplicative $\varepsilon^{-1}$ factor. The result follows from \cref{thm:Bound_EOT1}.
\begin{remark}[Extensions beyond
shallow ReLU]
    \begin{enumerate}
        \item[(1)] (Sigmoidal NNs). The results of this section readily extend to cover sigmoidal NNs, with slightly modified parameters. Specifically, one has to replace $s$ from \cref{thm:Bound_EOT1} with $\tilde s=\lfloor d_x / 2\rfloor+2$ and consider the sigmoidal NN class, with nonlinearity $\psi(z)=\left(1+e^{-z}\right)^{-1}$ (instead of ReLU) and parameters satisfying
\[
\max _{1 \leq i \leq k}\left\|w_i\right\|_1 \vee\left|b_i\right| \leq k^{\frac{1}{2}}\log k,\ \max _{1 \leq i \leq k}\left|\beta_i\right| \leq 2ak^{-1},\ \left|b_0\right| \leq a,\ \left\|w_0\right\|_1 = 0.
\]
The proofs of Theorems \ref{thm:Bound_EOT1}-\ref{thm:Bound_NeuralCoupling_EOT} go through using the second part of Proposition 10 from \cite{sreekumar2022neural}, which relies on controlling the so-called Barron coefficient (cf. \cite{barron1992neural,barron1993universal,yukich1995sup}). 
\item[(2)] (Deep NNs). Our theory currently accounts for NEs realized by shallow NNs, but deep nets are oftentimes preferable in practice. Extending our results to deep NNs should be possible by utilizing existing function approximation error bounds \cite{shen2019deep,schmidt2020nonparametric,bresler2020sharp}, although the resulting bounds do not seem to be sharp enough to yield the parametric rate of convergence. Nevertheless, we conjecture that deep NEs should achieve parametric convergence, and under milder regularity assumptions than their shallow counterparts. Sharper covering bounds would be needed to show that.

\item[(3)] (Kernel methods).  As an alternative to parametrizing the dual potential with NNs, one can use a function class based on a reproducing kernel Hilbert space (RKHS). This approach was proposed in \cite{nguyen2010estimating} for $f$-divergences, where an RKHS-based estimator was analyzed under the assumption that the optimizing dual function belongs to the RKHS, i.e., with no approximation error. This is a high-level assumption that is difficult to verify in practice, as no primitive sufficient conditions ensuring it were provided. The RKHS framework is particularly suitable for the EOT problem due to the intrinsic regularity of the EOT dual potentials. This regularity allows achieving the same performance guarantees as in \cite{nguyen2010estimating} while removing the aforementioned assumption. Specifically, \cite[Theorem 2]{genevay2019sample} shows that the optimal potential $\varphi$ of $\mathsf{OT}_c^{\varepsilon}(\mu, \nu)$, with $c \in C^\infty(\mathcal{X} \times \mathcal{Y})$, belongs to the Sobolev space $\mathbf{H}^s(\mathcal{X})$ and satisfies $\|\varphi\|_{\mathbf{H}^s} = O\!\left(1 + \varepsilon^{-(s-1)}\right)$, independently of the marginals. By the Sobolev embedding theorem, when $s > d_x / 2$, $\mathbf{H}^s(\mathcal{X})$ embeds continuously into an RKHS. We therefore define 
\[
\mathcal{G}_M := \{ h(x) : \|h\|_{\mathbf{H}^s} \le M \}, \quad M = O\!\left(1 + \varepsilon^{-(s-1)}\right),
\]
which by construction contains the true optimal potential. Replacing the NN class $\mathcal{F}_{k, d_x}(a)$ in \eqref{eq:NE_EOTc} with $\mathcal{G}_M$ eliminates the function approximation error, leaving only the statistical estimation error to be bounded. As the $L_2$ metric entropy of $\mathcal{G}_M$ scales as $O\!\left(\delta^{-d_x / s}\right)$ for $s > d_x / 2$~\cite{birman1967piecewise}, the same approximation–estimation decomposition as in the proof of \cref{thm:Bound_EOT1} with $s = \lceil d_x / 2 \rceil + 1$, yields the parametric estimation rate 
\[
O\!\left( \left(1 + \varepsilon^{-\lceil d_x / 2 \rceil}\right) n^{-1/2} \right).
\]
However, standard kernel estimators require manipulating Gram matrices, incurring $O(n^2)$ memory and $O(n^3)$ time complexity, which becomes prohibitive for large-scale datasets, and thus less scalable than NN parameterizations. Analogous arguments extend to the EGW estimation problem discussed in \cref{sec:Neural Estimation}.
    \end{enumerate}
    \end{remark}

\begin{remark}[Modular Analysis of Optimization Error]
\label{rem:opt_error}
    Our analysis in \cref{thm:Bound_EOT1} accounts for the function approximation and statistical estimation errors, but assumes an ideal optimizer in \eqref{eq:NE_EOTc}. While \emph{global convergence} guarantees for optimization algorithms of general NNs over nonconvex loss landscapes remains a challenging open problem, our analysis is modular enough to allow factoring out the optimization error as an additional additive term. Specifically, consider an arbitrary (possibly stochastic) optimization algorithm for solving the
problem in (13) from the datasets $(X^n,Y^n)$. This algorithm learns the NN parameters $\boldsymbol{\theta}_i:=(w_i,b_i,\beta_i)$. Following standard learning-theoretic formalism, we model such an algorithm by a transition kernel $P_{\boldsymbol{\theta}|X^n, Y^n}$. Write $\widehat{\mathsf{OT}}_{c,k,a}^{\mspace{1mu}\varepsilon}(X^n, Y^n; \boldsymbol{\theta})$ for the EOT value realized by the specific parameters $\boldsymbol{\theta}$ returned by the algorithm, and suppose that the global optimization error is bounded in expectation as
\begin{align*}
 \mathbb{E}_{X^n, Y^n}\left[\mathbb{E}_{\boldsymbol{\theta}|X^n, Y^n}\left[\left| \widehat{\mathsf{OT}}_{c,k,a}^{\mspace{1mu}\varepsilon}(X^n, Y^n; \boldsymbol{\theta}) - \widehat{\mathsf{OT}}_{c,k,a}^{\mspace{1mu}\varepsilon}(X^n, Y^n) \right|\right]\right] \le \delta_{\mathrm{opt}},   
\end{align*}
where $\widehat{\mathsf{OT}}_{c,k,a}^{\mspace{1mu}\varepsilon}(X^n, Y^n)$ is defined in \eqref{eq:NE_EOTc}.
Given the above, we can control the total error between the estimate obtained from the optimization algorithm $P_{\boldsymbol{\theta}|X^n, Y^n}$ and the ground-truth EOT value $\mathsf{OT}_c^{\varepsilon}(\mu,\nu)$ using the triangle inequality:
\begin{align*}
    &\mathbb{E}\left[\left| \widehat{\mathsf{OT}}_{c,k,a}^{\mspace{1mu}\varepsilon}(X^n, Y^n; \boldsymbol{\theta}) - \mathsf{OT}_c^{\varepsilon}(\mu,\nu) \right|\right] \\
    &\le \mathbb{E}\left[\left| \widehat{\mathsf{OT}}_{c,k,a}^{\mspace{1mu}\varepsilon}(X^n, Y^n; \boldsymbol{\theta}) - \widehat{\mathsf{OT}}_{c,k,a}^{\mspace{1mu}\varepsilon}(X^n, Y^n) \right|\right] + \mathbb{E}\left[\left| \widehat{\mathsf{OT}}_{c,k,a}^{\mspace{1mu}\varepsilon}(X^n, Y^n) - \mathsf{OT}_c^{\varepsilon}(\mu,\nu) \right|\right] \\
    &\le \delta_{\mathrm{opt}} + \delta_{\mathrm{EOT}},
\end{align*}
where the first term is bounded using the law of total expectation, while the second term denotes $\delta_{\mathrm{EOT}}$ the right-hand side (RHS) bound of \eqref{Bound_EOT1} in \cref{thm:Bound_EOT1}. Analogous arguments extend to the EGW estimation problem discussed in \cref{sec:Neural Estimation}.
\end{remark}

\section{Neural Estimation of the EGW Alignment Cost and Plan}\label{sec:Neural Estimation}

Building on the results of the previous section, we treat neural estimation of the EGW alignment cost and plan. Like before, we provide sharp bounds on the function approximation and statistical estimator errors, assuming an ideal optimizer for learning the neural dual potentials from the parametrized empirical objective. Without loss of generality, further suppose that $\mu$ and $\nu$ are centered (due to translation invariance) and that $d_x\leq d_y$. We start by describing the NEs for the EGW cost and plan, followed by non-asymptotic error bounds to guarantee their performance.  All proofs are deferred to the \cref{appen:Proofs}.

\subsection{EGW Neural Estimator}
\medskip
Consider the variational representation of the EGW distance stated in \cref{lemma:EGWDual}. The constant $\sC_{\mu,\nu}$ involves only moments of the marginal distributions and, as such, we estimate it using empirical averages via $\widehat{\sC}(X^n, Y^n)\coloneqq \sC_{\hat\mu_n,\hat\nu_n}$. For the second term in \eqref{eq:S2Decomp}, we estimate $\OT_\bA ^{\mspace{1mu}\eps}(\mu,\nu)$ using the EOT NE from \eqref{eq:NE_EOTc} and optimize over the auxiliary $\bA$ matrix using gradient descent. Overall, the resulting EGW estimator is given by 
\begin{equation}
\label{eq:NE_EGW}
\widehat{\mathsf{GW}}_{k,a}^{\mspace{1mu}\eps}(X^n, Y^n)\coloneqq\widehat{\sC}(X^n, Y^n)+\inf _{\bA  \in \cD_{M }}\Big\{32\|\bA \|_\F^2+\widehat\OT_{\bA , k, a}^{\mspace{1mu}\eps}(X^n, Y^n)\Big\}.
\end{equation}
where $M =\sqrt{d_xd_y}$. This choice is based on the fact that $\sqrt{d_xd_y} \geq M_{\mu, \nu} \vee M_{\hat\mu_n, \hat\nu_n}$ and the RHS of \eqref{eq:S2Decomp} achieves its global minimum inside $\mathcal{D}_{M_{\mu, \nu}}$. This estimator can be computed via alternating optimization over the matrices $\bA $ and the NN~parameters.

Akin to \eqref{eq:NE_EOTc_plan}, for any $(\bA ,f)\in\cD_{M }\times \cF_{k,d_x}(a)$, we define the induced neural alignment plan
\begin{equation}
    d \pi^\eps_{\bA ,f}(x, y) \coloneqq \frac{\exp \left(\frac{f(x)-c_\bA (x, y)}{\eps}\right)}{\int_{\cX} \exp \left(\frac{f(x)-c_\bA (x, y)}{\eps}\right) d\mu(x)} d\mu\otimes\nu(x,y).\label{eq:NE_plan}
\end{equation}
and write $\pi^\eps_{\bA_\star ,f_\star}$ for the plan induced by an optimal pair $(\bA _\star,f_\star)\in\cD_{M }\times \cF_{k,d_x}(a)$. The latter, provides a proxy of the true optimal plan $\pi^\eps_\star$ that achieves the infimum in \eqref{eq:QuadEGW}.

\subsection{Performance Guarantees}

Building on the error analysis for the NE of EOT with a general cost, we provide sharp non-asymptotic error bounds for the EGW problem. Our guarantees account for estimation of both the alignment cost and plan. 

\medskip
Starting from the EGW cost estimation question, we again establish two bounds, each with its utility. The first, presented below, holds for any dataset and NN sizes ($n$ and $k$, respectively) but contains a complicated dimension-dependent constant. We follow this result up with a second bound that rids of that constant, at the expense of an extra $\text{polylog}(k)$ factor and the requirement that $k$ is large enough.

\begin{theorem}[EGW cost neural estimation; bound 1] \label{thm:Bound_EGW1}
     There exists a constant $C>0$ depending only on $d_x,d_y$, such that setting $a=C(1+\eps^{1-s})$ with $s=\left\lfloor d_x  / 2\right\rfloor+3$, we have
\begin{equation}
\label{Bound_EGW1}
\begin{aligned}
&\sup _{(\mu, \nu) \in \mathcal{P}(\cX) \times \mathcal{P}(\cY)} \mathbb{E}\left[\left|\widehat{\mathsf{GW}}_{k,a}^{\mspace{1mu}\eps}(X^n, Y^n)-\mathsf{GW}^{\eps}(\mu, \nu)\right|\right]\\
&\qquad\qquad\lesssim_{d_x,d_y} \left(1+\frac{1}{\eps^{\left\lfloor \frac{d_x}{2}\right\rfloor+2}}\right) k^{-\frac{1}{2}}+\min\left\{1+\frac{1}{\eps^{\left\lceil d_x+\frac{d_y}{2}\right\rceil+4}}\,,\left(1+\frac{1}{\eps^{\left\lfloor \frac{d_x}{2}\right\rfloor+2}}\right)\sqrt{k}\right\} n^{-\frac{1}{2}}.
\end{aligned}
\end{equation}
\end{theorem}
The proof of \cref{thm:Bound_EGW1} is given in \cref{appen:proof_Bound_EGW1}. Thanks to the EGW variational form given in \cref{lemma:EGWDual} and the structure of the NE, the error analysis decomposes as 
\begin{align*}
&\mathbb{E}\left[\left|\widehat{\mathsf{GW}}_{k,a}^{\mspace{1mu}\eps}(X^n, Y^n)-\mathsf{GW}^{\eps}(\mu, \nu)\right|\right]\\
&\qquad\qquad\qquad\leq \mathbb{E}\big[\big|\mathsf{C}_{\mu, \nu}-\widehat{\sC}(X^n, Y^n)\big|\big]+\mathbb{E}\left[\sup_{\bA\in\cD_M}\left|\OT^\eps_{c_{\bA} ,k,a}(\mu, \nu)-\widehat\OT^\eps_{c_{\bA} , k,a}\left(X^n, Y^n\right)\right|\right],
\end{align*}
where error for $\mathsf{C}_{\mu, \nu}$ is easy to analyze (namely, estimation of certain marginal moments via sample means), while the second term requires controlling the estimation error of the EOT cost $\OT_{\bA}^{\mspace{1mu}\eps}(\mu,\nu)$, uniformly in $\bA \in\cD_{M }$. The latter decomposed into a (worst-case in $\bA$) approximation and estimation error terms, as in \eqref{eq:EOT_NE_error_decomposition}. Each is then bounded by invoking the results of Lemmas \ref{lem:approximation_error_bound} and \ref{lem:estimation_error_bound} with the cost function $c_\bA$, and controlling cost-dependent terms uniformly in $\bA\in\cD_M$, e.g., in terms of properties of the marginal populations $\mu$ and $\nu$. As in the EOT case, the resulting bound requires calibrating the parameter $a$, which controls the magnitude of the NN weights, to the dimension-dependent constant $C$ (see Remarks~\ref{REM:C_constant}~and~\ref{REM:C_constant_egw}).

\begin{remark}[Almost explicit expression for $C$ and dependence on dimension]
\label{REM:C_constant_egw}
The expression of the constant $C$ in \cref{thm:Bound_EGW1} can also be evaluated as in Remark 3. We can write $C=C_sC_{ d_x, d_y} \bar{c}_{d_x}$, with explicit expressions for $C_{d_x, d_y}$ and $\bar{c}_{d_x}$ given in \eqref{eq:constant_C_dx_dy_egw} and \eqref{eq:constant_c_bar}, respectively, while $C_s$ is a combinatorial constant that arises from the multivariate Faa di Bruno formula (cf. \eqref{eq:constant_C_s_1}-\eqref{eq:constant_C_s_2}). The dependence of $C$ on the dimension, which is viewed as fixed in this work, is (super-)exponential. This is primarily because the smoothness index $s=s_{\mathsf{KB}} := \lfloor d_x/2 \rfloor + 3$ features in the exponent of $\mathrm{poly}(d_x)$ terms within both $\bar{c}_{d_x}$ and $C_{d_x, d_y}$, and also influences the combinatorial term $C_s$. This structure leads to a dependence of the form $d_x^{O(d_x)}$. The convoluted nature of $C_s$ is the main reason we view the overall constant $C$ as implicit.
\end{remark}

\begin{remark}[Optimality and dependence on dimension]
    By equating $k\asymp n$, the proposed EGW NE achieves the parametric $n^{-1/2}$ rate, and is thus minimax rate-optimal. Further observe that the bound can be instantiated to depend only on the smaller dimension $d_x$ (by omitting the first argument of the minimum). This is significant when one of the mm spaces has a much smaller dimension than the other, as the bound will adapt to the minimum. This phenomenon is called the LCA principle and it was previously observed for the empirical plug-in estimator of the OT cost \cite{hundrieser2022empirical}, the EOT and EGW distances \cite{groppe2023lower}, and the GW distance itself \cite{zhang2024gromov}.
\end{remark}

The second bound removes the dependence on $C$, analogously to \cref{thm:Bound_EOT_2}.

\begin{theorem}[EGW cost neural estimation; bound 2] \label{thm:Bound_EGW2}
    Let $\epsilon>0$ and set $m_k=\log k\vee 1$. Assuming $k$ is sufficiently large, we have
\begin{equation}
\label{Bound_EGW2}
\begin{aligned}
   &\sup _{(\mu, \nu) \in \mathcal{P}(\cX) \times \mathcal{P}(\cY)} \mspace{-5mu}\mathbb{E}\left[\left|\widehat{\mathsf{GW}}_{k,m_k}^{\mspace{1mu}\eps}(X^n, Y^n)\mspace{-2mu}-\mspace{-2mu}\mathsf{GW}^{\eps}(\mu, \nu)\right|\right]\\
   &\qquad\qquad\qquad\lesssim_{d_x,d_y} \left(1+\frac{1}{\eps^{\left\lfloor\frac{d_x}{2}\right\rfloor+2}}\right) k^{-\frac{1}{2}}+\min\left\{\left(1+\frac{1}{\eps^{\left[\frac{d_y}{2}\right\rceil}}\right)(\log k)^2\,, \sqrt{k}\log k\right\} n^{-\frac{1}{2}}.
\end{aligned}
\end{equation}
\end{theorem}

The proof is almost identical to that of \cref{thm:Bound_EOT_2}, and is omitted to avoid repetition. We note that a partial account of the requirement that $k$ is large enough can be provided. Specifically, for the bound to hold we need $k$ to be such that $\log k \geq C(1+\eps^{1-s})$, where $C$ is the constant from \cref{thm:Bound_EGW1}. However, it is challenging to provide a simple expression for the threshold on $k$ due to the implicit nature of $C$.

\medskip
Lastly, we account for the quality of the neural alignment plan from  \eqref{eq:NE_plan} by comparing it, in KL divergence, to the true EGW alignment $\pi^\star$.
\begin{theorem}[EGW alignment plan neural estimation]
\label{thm:Bound_NeuralCoupling}
    Suppose that $\mu \in \mathcal{P}_{\mathsf{ac}}(\cX)$. Let $\bA _\star$ be a minimizer of \eqref{eq:S2Decomp} and $\hat{f}_\star$ be a maximizer of $\widehat\OT^\eps_{\bA _\star,k,a}(X^n,Y^n)$ from \eqref{eq:NE_EGW}, with $a$ as defined in \cref{thm:Bound_EGW1}. Then, the induced neural alignment plan $\pi^\eps_{\bA _\star,\hat{f}_\star}$ from \eqref{eq:NE_plan} satisfies
    \begin{equation}
    \label{Bound_NeuralCoupling}
\mathbb{E}\left[\KL\left(\pi_\star^\eps\middle\|\pi^\eps_{\bA _\star,\hat{f}_\star}\right)\right] \lesssim_{d_x,d_y} \varepsilon^{-1}\delta_{\mathrm{EGW}}
    \end{equation}
where $\pi^\eps_\star$ is optimal coupling of EGW problem \eqref{eq:QuadEGW} and  $\delta_{\mathrm{EGW}}$ the right hand side bound of \eqref{Bound_EGW1} in \cref{thm:Bound_EGW1}.
\end{theorem}

The proof is almost identical to that of \cref{thm:Bound_NeuralCoupling_EOT} since we account for optimal $\bA_*$ here (i.e. $c=c_{\bA_\star}$) and Corollary 1 from \cite{rioux2024entropic} implies that $\pi^\eps_{\bA_\star} = \pi^\eps_\star$, where $\pi^\eps_{\bA_\star}$ is the unique EOT coupling for 
$\OT^\eps_{\bA _\star}(\mu,\nu)$, then the result follows by invoking \cref{thm:Bound_NeuralCoupling_EOT} and \cref{thm:Bound_EGW1}.

\begin{remark}[Limitation of \cref{thm:Bound_NeuralCoupling}]
    We note that \cref{thm:Bound_NeuralCoupling} does not fully account for the quality of learned (empirical) neural alignment plan since the statement considers an optimal $\mathbf{A} _\star$, which is unknown in practice. This limitation comes from our analysis, which can only account for couplings corresponding to the same matrix $\mathbf{A}$ on both sides of the KL divergence. One can show that the gap between the KL divergence of the empirically optimal plan (using $\hat{\mathbf{A}}_{\star}$) and the surrogate bounded by our theorem (using $\mathbf{A}_{\star}$) is controlled by:
\[
\mathbb{E}\mspace{-2mu}\left|\mathsf{D}_{\mathsf{KL}}\left(\mspace{-4mu}\pi_{\star}^{\varepsilon} \middle \| \pi_{\hat{\mathbf{A}}_{\star}, \hat{f}_{\hat{\mathbf{A}}_\star}}^{\varepsilon}\mspace{-4mu}\right)\mspace{-4mu}-\mspace{-2mu}\mathsf{D}_{\mathsf{KL}}\left(\pi_{\star}^{\varepsilon} \middle \| \pi_{\mathbf{A}_{\star}, \hat{f}_{\mathbf{A}_\star}}^{\varepsilon}\right)\right|\mspace{-2mu}\leq\mspace{-2mu} \frac{2}{\varepsilon}\mathbb{E}\mspace{-3mu}\left[\mspace{-2mu}\left(\mspace{-2mu}\left\|\hat{f}_{\mathbf{A}_\star}\mspace{-5mu}-\mspace{-3mu}\hat{f}_{\hat{\mathbf{A}}_\star}\right\|_{\infty, \mathcal{X}}\mspace{-9mu}+\mspace{-2mu}32 \sqrt{d_x d_y}\left\|\hat{\mathbf{A}}_\star\mspace{-5mu}-\mspace{-3mu}\mathbf{A}_\star\right\|_F\mspace{-2mu}\right)\mspace{-2mu}\right].
\]
This bound reduces the problem to proving the stability of the learned neural network potential with respect to $\mathbf{A}$. Strengthening our result to formally account for the empirically optimized matrix $\hat{\mathbf{A}}_{\star}$ is of significant interest. This, however, would require a delicate stability analysis of the learned potential within the EGW problem, which we leave as an important direction for future work.
\end{remark}


\section{Computation and Numerical Experiments}\label{sec:experiments}

Variants of neural EOT estimation were previously used in practice on different tasks, demonstrating scalability and accuracy \cite{seguy2018large,daniels2021score,mokrov2023energy}. The EGW NE, on the other hand, is novel and has not been empirically tested before. This section provides numerical experiments for EGW neural estimation over synthetic as well as real-world data. The results validate our theory and demonstrate similar virtues to its EOT counterpart. We start by describing the specifics of the algorithm and then present the results.

\subsection{Neural Estimation Algorithm}\label{subsec:experiments_algorithm}

\begin{algorithm}[!t]
\caption{Entropic Gromov-Wasserstein distances via Neural Estimation}\label{alg:EGW_NE}
\begin{algorithmic}[1]
\Statex \textbf{Input}: Entropic parameter $\eps >0$, epoch number $N > 0$, batch size $m >0$, learning rate $\eta >0$, initial matrix $\bm C_0\in\mathcal D_{M}$ for any $M \geq M_{\hat\mu_n,\hat\nu_n}$ and initialized neural network $\hat f_0$. Fix the step sequences $\beta_k=\frac{1}{2L}$, $\gamma_k=\frac{k}{4L}$, and $\tau_k=\frac{2}{k+2}$, with $L=64 \vee\left(32^2 \eps^{-1} \sqrt{M_4\left(\hat\mu_n\right) M_4\left(\hat\nu_n\right)}-64\right)$.
\State $k \gets 1$
\State $\bm A_1\gets \bm C_0$
\State $\hat{f}_1 \gets \text{Adam}(\hat f_0,\bA _1;N,m,\eta)$
\State $\bm G_1 \gets 64 \bA _1-32 \sum_{1 \leq i,j \leq n} X_iY_j^{\intercal} \big[\widetilde{\bm{\Pi}}^\eps_{\bA _1}\big]_{ij}$, 
\While{stopping condition is not met}
\State $\bm B_{k}\gets \frac M2 \sign( \bm A_k-\beta_k\bm G_k)\min\left(\frac{2}{M}\left|\bm A_k-\beta_k\bm G_k\right|,1\right)$
\State $\bm C_{k}\gets \frac{M}{2}\sign(\bm C_{k-1}-\gamma_k\bm G_k)\min\left(\frac{2}{M}\left|\bm C_{k-1}-\gamma_k\bm G_k\right|,1\right)$
\State $\bm A_{k+1}\gets\tau_k\bm C_k+(1-\tau_k) \bm B_k$
\State $\hat{f}_{k+1} \gets \text{Adam}(\hat f_k,\bA _{k+1};N,m,\eta)$
\State $\bm G_{k+1} \gets 64 \bA _{k+1}-32 \sum_{1 \leq i,j \leq n} X_iY_j^{\intercal} \big[\widetilde{\bm{\Pi}}^\eps_{\bA _{k+1}}\big]_{ij}$, 
\State $k\gets k+1$
\EndWhile
\State \textbf{Output}: $(\bm B_k,\hat{f}_{k+1})$

\end{algorithmic}
\end{algorithm}

Our goal is to compute $\widehat{\mathsf{GW}}_{k,a}^{\mspace{1mu}\eps}(X^n, Y^n)$ from \eqref{eq:NE_EGW}. The constant $\widehat{\sC}(X^n,Y^n)$ is computed offline, so we focus  on solving 
\begin{equation}
\label{eq:alg_objective}
\inf _{\bA  \in \cD_{M }}32\|\bA \|_\F^2+\sup_{f \in \cF_{k, d_x}(a)} \left\{ \frac{1}{n}\sum_{i=1}^nf(X_i)-\frac{\eps}{n}\sum_{j=1}^n\log\left(\frac{1}{n}\sum_{i=1}^n\exp\left(\frac{f(X_i)-c_{\bA}(X_i,Y_j)}{\eps}\right)\right)\right\},
\end{equation}
where $M =\sqrt{d_xd_y}$ (which guarantees that all the optimizers are within the optimization domain; cf. \cite[Corollary 1]{rioux2024entropic}) and $a$ is unrestricted, so as to enable optimization over the whole parameter space. 

For the outer minimization, we employ the accelerated first-order method with an inexact oracle from \cite{rioux2024entropic}. Algorithms 1 therein accounts for the case when $\eps>16\sqrt{M_4(\mu)M_4(\nu)}$, whence the optimization objective becomes convex in $\bA $ (and is always smooth). Theorem 9 of that work guarantees that this algorithm converges to a global solution up to an error that is linear in the precision parameter of the inexact oracle. The convergence follows the optimal $O(k^{-2})$ rate for smooth constrained optimization \cite{nesterov2003introductory}, with $k$ being the iteration count. Algorithm 2 of \cite{rioux2024entropic} relies on similar ideas but does not assume convexity, at the cost of local convergence guarantees and a slower rate of $O(k^{-1})$. In our experiments, we find that the convexity condition on $\eps$ is quite restrictive, and hence mostly employ the latter method.

Each iteration of Algorithm 2 from \cite{rioux2024entropic}, requires computing the approximate gradient 
\begin{equation}
    64 \bA -32 \sum_{1 \leq i,j \leq n} X_iY_j^{\intercal} \big[\widetilde{\bm{\Pi}}^\eps_{\bA}\big]_{ij},\label{eq:EGW_approx_grad}
\end{equation} 
where $\widetilde{\bm{\Pi}}^\eps_{\bA}\in\RR^{n\times n}$ is the oracle coupling (since we are in the discrete setting, the coupling is given by a matrix and is denoted as such). We instantiate the oracle as the neural estimator, corresponding to the inner optimization in \eqref{eq:alg_objective}. Specifically, for a fixed $\bA $, we train the parameters of the ReLU network using the Adam algorithm \cite{kingma2014adam}. Upon obtaining an optimizer $\hat{f}_\bA $, we compute the induced neural coupling
\begin{equation}
\label{eq:empirical_neural_coupling}\big[\widetilde{\bm{\Pi}}^\eps_{\bA}\big]_{ij}=\frac{\exp \left(\frac{\hat{f}_\bA \left(X_i\right)-c_{\bA}\left(X_i, Y_j\right)}{\eps}\right)}{n \sum_{k=1}^n \exp \left(\frac{\hat{f}_\bA \left(X_k\right)-c_{\bA}\left(X_k, Y_j\right)}{\eps}\right)}, \quad i,j=1,\ldots,n.   
\end{equation}
This coupling is plugged into the aforementioned approximate gradient formula, and we alternate between the optimization over the matrices $\bA $ and the NN $f$ until convergence. The method is summarized in \cref{alg:EGW_NE}.

The convergence guarantees from \cite{rioux2024entropic} (specifically, their Theorem 11 on adaptive convergence rate and Corollary 12 on approximate stationarity) apply directly to our \cref{alg:EGW_NE}, provided that our neural coupling $\widetilde{\boldsymbol{\Pi}}^\varepsilon_{\mathbf{A}}$ from \eqref{eq:empirical_neural_coupling} is a $\delta$-oracle approximation of the true empirical coupling $\boldsymbol{\Pi}^\varepsilon_{\mathbf{A}}$ of $\mathsf{OT}_{\mathbf{A}}^\varepsilon\left(X^n, Y^n\right)$. The following lemma shows that if the NN obtained from optimizing \eqref{eq:alg_objective}  is close to an optimal dual potential of $\mathsf{OT}_{\mathbf{A}}^\varepsilon\left(X^n, Y^n\right)$, then the resulting neural coupling is a $\delta$-oracle.

\begin{lemma}
\label{lem:delta-oracle}
    Suppose that for any fixed $\mathbf{A}$,  the inner maximization in \eqref{eq:alg_objective}  achieves global optimality up to error $\delta^\prime$, i.e., $\|\hat{f}_{\mathbf{A}}-\varphi_{\mathbf{A}}\|_{\infty,\mathcal{X}} \leq \delta^\prime $, where $\hat{f}_{\mathbf{A}}$ is the NN returned by our algorithm solving the inner maximization and $\varphi_{\mathbf{A}}$  is the optimal dual potential. Then, 
    \begin{equation}
        \label{eq:delta-oracle}
        \left\|\widetilde{\boldsymbol{\Pi}}^\varepsilon_{\mathbf{A}}-\boldsymbol{\Pi}^\varepsilon_{\mathbf{A}}\right\|_{\infty}\leq \frac{\delta^\prime}{2n \varepsilon}:=\delta.
    \end{equation}
\end{lemma}
The proof is given in \cref{appen:proof:delta-oracle}. \cref{lem:delta-oracle} reduces the problem of constructing accurate proxies for couplings to constructing accurate dual potentials for any fixed $\mathbf{A}$. To fill this gap, one need to etablish quantitative stability bounds of dual potentials in $\mathbf{A}$, which is a hard problem that remains open.

        

\begin{remark}[Computational and memory complexity]
\label{rmk:complexity}
    The proposed mini-batch NE achieves substantial computational and memory gains over standard entropic OT solvers such as Sinkhorn's algorithm. Let $b$ denote the mini-batch size and $k$ the number of neurons in the shallow network. Each gradient step involves two main operations:
\begin{itemize}
    \item[(i)] OT objective evaluation: forming the $b\times b$ matrix of pairwise interactions for the log-sum-exp term, which costs $O(b^2)$ time;
    \item [(ii)] NN pass: forward and backward propagation through the network, which costs $O(bk)$ time.
\end{itemize}
Hence, the per-step cost is $O(b^2+bk)$ and the per-epoch cost scales as
\[
O\!\left(\frac{n}{b}(b^2+bk)\right)=O\!\left(n(b+k)\right).
\]
In the large-scale regime, both $b$ and $k$ are fixed hyperparameters independent of $n$, yielding overall linear-time complexity $O(n)$ per epoch and memory requirement $O(b^2+bk)$. 
In contrast, the classical Sinkhorn algorithm \cite{cuturi2013sinkhorn} must form and iterate over the full $n\times n$ cost matrix, incurring $O(n^2)$ time and memory per iteration. This quadratic scaling renders Sinkhorn infeasible for large datasets, while our NE scales efficiently and maintains stable empirical performance even for high-dimensional and large-sample settings.
\end{remark}

\subsection{Synthetic Data}

\begin{figure}[t!]
  \begin{subfigure}{0.345\textwidth}
    \centering
    \includegraphics[width=\textwidth]{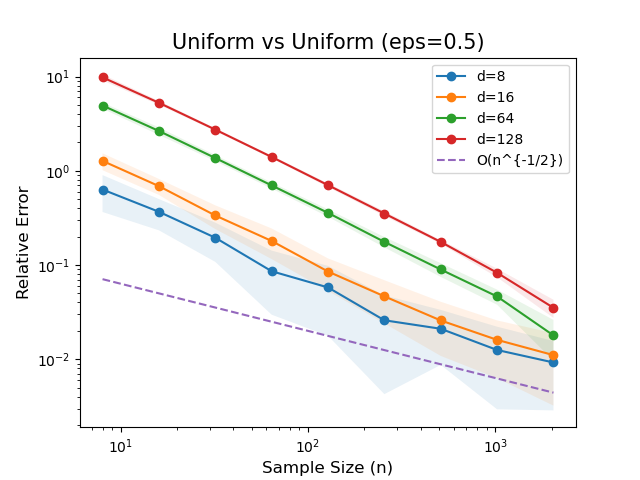}
    \caption{}
    \label{fig:uniform}
  \end{subfigure}%
  \hspace{-6mm}
  \begin{subfigure}{0.345\textwidth}
    \centering
    \includegraphics[width=\textwidth]{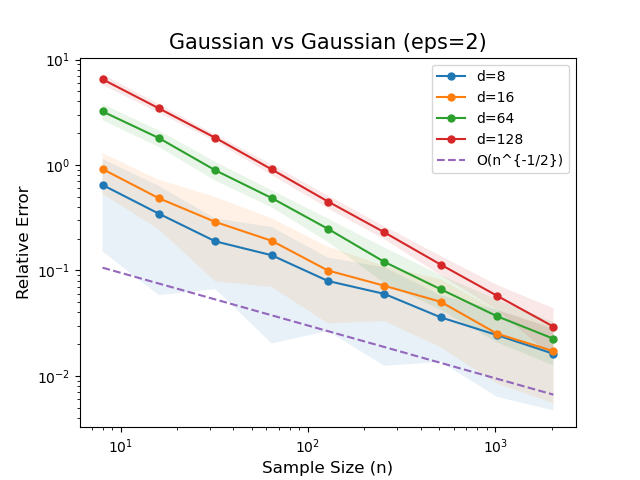}
    \caption{}
    \label{fig:gaussian}
  \end{subfigure}
  \hspace{-6mm}
    \begin{subfigure}{0.345\textwidth}
    \centering
    \includegraphics[width=\textwidth]{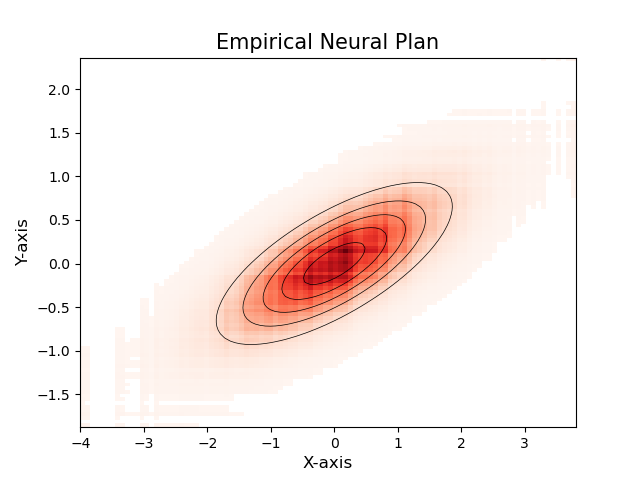}
    \caption{}
    \label{fig:Neural_Plan}
  \end{subfigure}%
  \caption{Neural Estimation of EGW alignment: (a) Relative error for the case where $\mu=\nu=\mathrm{Unif}\big([-1/\sqrt{d},1/\sqrt{d}]^d\big)$; (b) Relative error for $\mu,\nu$ as centered Gaussian distributions with randomly generated covariance matrices; (c) Learned neural alignment plan (in red) versus the true optimal GW alignment (whose density is represented by the back contour lines).}
  \label{fig:simulated_data}
\end{figure}

We test our neural estimator on synthetic data, by estimating the EGW cost and alignment plan between uniform and Gaussian distribution in different dimensions. We consider dimensions $d\in \{8,16,64,128\}$, and for each $d$, employ a ReLU network of size $k\in \{32,64,128,256\}$, respectively. Accuracy is measured using the relative error $\big|\widehat{\mathsf{GW}}_{k,a}^{\mspace{1mu}\eps}(X^n, Y^n)-\widetilde{\mathsf{GW}}^\eps(\mu,\nu)\big|/\,\widetilde{\mathsf{GW}}^\eps(\mu,\nu)$, where $\widetilde{\mathsf{GW}}^\eps(\mu,\nu)$ is regarded as the group truth, which we obtain by running  \cite[Algorithm 2]{scetbon2022linear} with $n=10,000$ samples (which we treat as $n\to\infty$ as it is $\times 5$ more than the largest sample set we use for our neural estimator).\footnote{We benchmark against the algorithm from \cite{scetbon2022linear} due to its efficient memory usage. We have also attempted approximating the ground truth using Algorithm 2 from \cite{rioux2024entropic}, which relies on the Sinkhorn oracle, but could not obtain stable results when scaling up to larger $d,n$ values.} Each of the presented plots is averaged over 20 runs. 

We first consider the EGW distance with $\eps=0.5$ between two uniform distribution over a hypercube, namely,  $\mu=\nu=\mathrm{Unif}\big([-1/\sqrt{d},1/\sqrt{d}]^d\big)$. \cref{fig:uniform} plots the EGW neural estimation error versus the sample size $n\in\{8,16,32,64,128,256,512,1024,2048\}$ in a log-log scale. The curves exhibit a slope of approximately $-1/2$ for all dimensions, conforming with our theory. In this experiment, we use 5 training epochs, and set the stopping condition for updating $\bA $ as either reaching a maximal iteration count of 100 or the Frobenius norm of gradient approximation dropping below $10^{-4}$. 

Next, we test the EGW NE on unbounded measures. To that end, we set $\eps = 2$ and take $\mu,\nu$ as centered $d$-dimensional Gaussian distributions with randomly generated covariance matrices. Specifically, the two covariance matrices are of the form $\mathbf{B}^{\intercal}\mathbf{B}+1/(3d)\mathbf{I}_d$, where $\mathbf{I}_d$ is a $d\times d$ identity and $\mathbf{B}$ is a matrix whose entries are randomly sampled from $\mathrm{Unif}([-1/d,1/d])$. Note that the generated covariance matrix is positive semi-definite with eigenvalues set to lie in $[\frac{1}{3d},\frac{1}{d}]$. \cref{fig:gaussian} plots the relative EGW neural estimation error for this Gaussian setting, again showing a parametric convergence rate. This experiment uses 10 epochs for the NE training, and the stopping condition for updating $\bA $ is either reaching 200 iterations or that the Frobenius norm of gradient approximation drops below $10^{-3}$. 

Lastly, we assess the quality of the neural alignment plan learned from our NE. Since doing so requires knowledge of the true (population) alignment plan $\pi^\eps_\star$, we consider the EIGW distance $\mathsf{IGW}^{\eps}(\mu,\nu)$ (i.e., with inner product cost; see \eqref{eq:EIGW}) between Gaussians, for which a closed form expression for the optimal plan was derived in \cite{le2022entropic}. Adapting \cref{alg:EGW_NE} (which treats EGW with quadratic cost) to the EIGW case merely amounts to changing the constants in the approximate gradient formula from the expression in \eqref{eq:EGW_approx_grad} to 
\[16 \bA -8 \sum_{1 \leq i,j \leq n} X_iY_j^{\intercal} \big[\widetilde{\bm{\Pi}}^\eps_{\bA}\big]_{ij},\]
where $\widetilde{\bm{\Pi}}^\eps_{\bA}$ is now the optimal coupling for the EOT problem with the cost function $(x,y)\mapsto -8x^{\intercal}\bA  y$ (see \cref{appen:EIGW} for more details on the EIGW setting). We take $\eps=0.5$, $\mu=\mathcal{N}\left(0,1\right)$, and $\nu=\mathcal{N}(0,1/4)$. By Theorem 3.1 of \cite{le2022entropic}, we have that $\pi^\eps_\star=\mathcal{N}\left(\mathbf{0}, \Sigma_\star\right)$ with
\[
\Sigma_\star=\left(\begin{array}{cc}
1 & \frac{1}{\sqrt{8}} \\
\frac{1}{\sqrt{8}} & \frac{1}{4}
\end{array}\right)
\]
Figure \ref{fig:Neural_Plan} compares the neural coupling learned from our algorithm, shown in red, to the optimal $\pi^\eps_\star$ given above, whose density is represented by the black contour line. The neural coupling is learned using $n=10^4$ samples and is realized by a NN with $k=40$ neurons. There is a clear correspondence between the two, which supports the result of~\cref{thm:Bound_NeuralCoupling}.\footnote{While \cref{thm:Bound_NeuralCoupling} is stated for quadratic EGW problem, the same conclusion holds true under the EIGW setting; see \cref{appen:EIGW}.} 

\subsection{MNIST Dataset}

\begin{figure}[t!]
  \begin{subfigure}{0.5\textwidth}
    \centering
    \includegraphics[width=\textwidth]{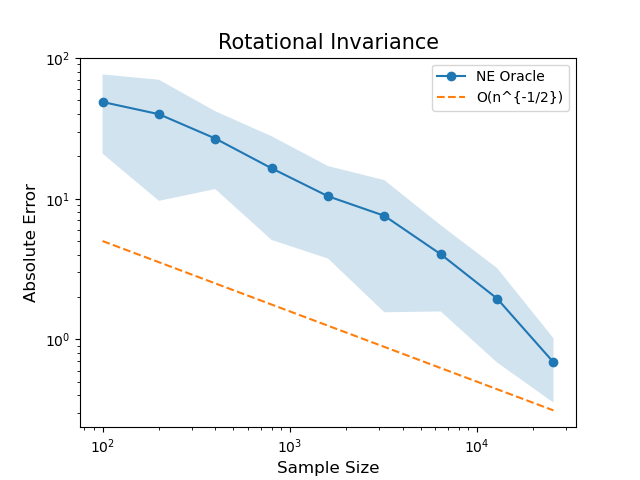}
    \caption{}
    \label{fig:NE_vs_Sinkhorn}
  \end{subfigure}%
  \begin{subfigure}{0.5\textwidth}
    \centering
    \includegraphics[width=\textwidth]{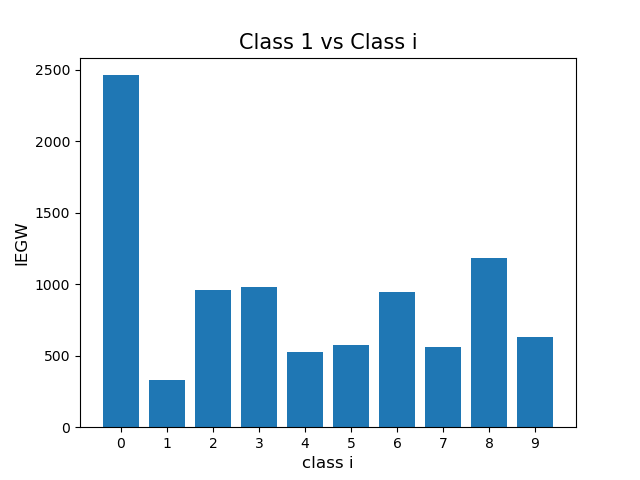}
    \caption{}
    \label{fig:Class_1_vs_i}
  \end{subfigure}
  \caption{Neural Estimation of EIGW on MNIST: (a) Testing for the orthogonal invariance of the EIGW distance by estimating the gap $\big|\mathsf{IGW}^\eps(\mu,\mu)-\mathsf{IGW}^\eps(\mathbf{U}_{\sharp}\mu,\mathbf{V}_{\sharp}\mu)\big|$, for $\mu$ as the empirical MNIST distribution and $(\mathbf{U},\mathbf{V})$ two orthogonal matrices; (b) Capturing visual similarities between digits by estimating the EIGW distance between different MNIST digits.}
  \label{fig:mnist_experiment}
\end{figure}

We next test our NE on the MNIST dataset, as a simple example of real-world data. We again consider the EIGW distance $\mathsf{IGW}^{\eps}(\mu,\nu)$ for these experiments, due to the improved numerical stability it provides. We also set $\eps$ to be large enough so that the algorithms does not incur numerical errors. To that end, we initiate a small $\eps$ value, and if errors occur, double it until the algorithm converges without
errors (eventually, we ended up using $\eps=10^3$, which is an order of magnitude smaller than the threshold for convexity condition from \cite[Theorem 2]{rioux2024entropic} to hold). 

We consider two experiments under the MNIST setting, one quantitative and another qualitative. For the first, we numerically test for the rotation invariance of the EIGW distance. Denoting the empirical distribution of the MNIST dataset by $\mu\in\cP(\RR^{784})$, for any two orthogonal matrices $\mathbf{U},\mathbf{V}$, we have $\mathsf{IGW}^{\eps}(\mu,\nu)=\mathsf{IGW}^{\eps}(\mathbf{U}_{\sharp}\mu,\mathbf{V}_{\sharp}\mu)$. Thus, we estimate the difference between these two quantities and evaluate its gap from 0, which is the ground truth. Specifically, we randomly sample data $(X_0^n,X_1^n)$ from $\mu\otimes\mu$, and $(Y_0^n,Y_1^n)$ from $(\mathbf{U},\mathbf{V})_{\sharp}\mu\otimes\mu$, with sample size $n\in\{100,200,400,800,1600,3200,6400,12800,25600\}$, and compute $\big|\mathsf{IGW}^\eps_{k,a}(X_0^n,X_1^n)-\mathsf{IGW}^\eps_{k,a}(Y_0^n,Y_1^n)\big|$ with $k=256$. 

\cref{fig:mnist_experiment} plots the estimation error (averaged over 10 runs), decreasing from 48 to 0.68, which again exhibits a parametric rate of convergence, as expected. At $n=25600$, the absolute error of about $0.6$ amounts to $10^{-3}$ relative error (compared to the value of the NE, which is around 512). Notably, despite the data dimension being relatively large in this experiment, the NE successfully retrieves the correct value with moderate sample sizes. We have also attempted replicating this experiment using the Sinkhorn-based oracle from \cite[Algorithm 2]{rioux2024entropic}, but the run time was prohibitively long on our devices for the larger $n$ values considered.

Lastly, we examine the ability of the EIGW NE to capture qualitative visual similarities between different digits in the MNIST dataset. To that end, we perform neural estimation of the EIGW distance between Class 1 and the other MNIST Classes, expecting, for example, a greater similarity (viz. smaller EIGW value) with digits like 4 and 7, but less so with 0, 3, and 8. The obtained estimates are plotted in \cref{fig:Class_1_vs_i}, where the expected qualitative behavior is indeed observed.

\section{Proofs}
\label{appen:Proofs}

We first introduce a technical result from approximation theory that will be used in the subsequent derivations. The following result, which is a restatement of Proposition 10 from \cite{sreekumar2022neural}, states that a sufficiently smooth function over a compact domain can be approximated to within $O(k^{-1/2})$ error by a shallow NN.

\begin{proposition}[Approximation of smooth functions; Proposition 10 from \cite{sreekumar2022neural}]
\label{prop:approx_restated}
Let $\cX\subseteq\RR^d$ be compact and $g:\cX \rightarrow \RR$. Suppose that there exists an open set $\cU \supset  \cX$, $b \geq 0$, and $ \tilde g \in  \cC_b^{s_{\mathsf{KB}}}(\cU)$, $s_{\mathsf{KB}}\coloneqq \lfloor d/2\rfloor+3$, such that $g= \tilde g\big|_\cX$.  
Then, there exists $f \in \cF_{k,d}\left(\bar c_{b,d,\|\cX\|}\right)$, where $\bar c_{b,d,\|\cX\|}$ is given in Equation (A.15) of \cite{sreekumar2022neural}, such~that 
\[\|f-g\|_\infty\lesssim \bar{c}_{b,d,\|\cX\|}d^{\frac 12} k^{-\frac 12}.\] 
\end{proposition}  

This proposition will allow us to control the approximation error of the NE. To invoke it, we will establish smoothness of the semi-dual EOT potentials arising in the variational representation from \eqref{eq:S2Decomp} (see \cref{lemma:regularity_EOT} ahead). The smoothness of potentials stems from the presence of the entropic penalty and the smoothness of the cost function~$c$. 

\subsection{Proof of Lemma \ref{lem:approximation_error_bound}}
\label{appen:proof_approx_error_bound}

\medskip
Proposition \ref{prop:approx_restated} provides a sup-norm approximation error bound of a smooth function by a NN. To invoke it, we first study the semi-dual EOT potentials and show that they are indeed smooth functions, i.e., admit an extension to an open set with sufficiently many bounded derivatives. The following lemma establishes regularity of semi-dual potentials for $\OT_c^{\eps}(\mu,\nu)$; after stating it we shall account for the extension. 

\begin{lemma}[Uniform regularity of EOT potentials]\label{lemma:regularity_EOT}
There exist semi-dual EOT potentials $(\varphi,\varphi^{c,\eps})$ for $\OT_c^{\mspace{1mu}\eps}(\mu, \nu)$, such that
\begin{equation}
\begin{aligned}
\label{eq:regularity_EOT}
    \|\varphi\|_{\infty,\cX} &\leq \frac{5}{2}\|c\|_{\infty, \mathcal{X} \times \mathcal{Y}}\\
    \left\|D^\alpha \varphi\right\|_{\infty,\cX} &\leq C_{s} C_{c,\|\mathcal{X}\|,\|\mathcal{Y}\|}^{s}\left(1+\eps^{1-s}\right) \text { with } 1 \leq|\alpha| \leq s,
\end{aligned}
\end{equation}
for any $s \geq 2$ and some constant $C_s$ that depends only on $s$ and $C_{c,\|\mathcal{X}\|,\|\mathcal{Y}\|}$ that depends only on cost $c(x,y)$ and $\|\mathcal{X}\|,\|\mathcal{Y}\|$. Analogous bounds hold for $\varphi^{c,\eps}$.
\end{lemma}
The lemma is proven in Appendix \ref{appen:proof:lemma:regularity_EOT}. The derivation is similar to that of Lemma 4 in \cite{zhang2024gromov}, but the bounds are adapted to the compactly supported case and present an explicit dependence on $\eps$ (as opposed to the $\eps=1$ assumption that was imposed in that work). 

Let $(\varphi,\varphi^{c,\eps})$ be semi-dual potentials as in  \cref{lemma:regularity_EOT} (i.e., satisfying \eqref{eq:regularity_EOT}) with the normalization $\int \varphi d \mu=\int \varphi^{c,\eps} d \nu=\frac{1}{2} \OT_c^{\eps}(\mu, \nu)$. Define the natural extension of $\varphi$ to the open ball of radius $\sqrt{2d_x}\|\mathcal{X}\|$:
\[
\tilde\varphi(x)\coloneqq -\eps \log \int_{\cY} \exp \left(\frac{\psi(y)-c(x, y)}{\eps}\right) d \nu(y), \quad x \in B_{d_x}\left(\sqrt{2d_x}\|\mathcal{X}\|\right),
\]
and notice that $\tilde\varphi\big|_\cX=\varphi$, pointwise on $\cX$. Similarly, consider its $(c,\eps)$-transform $\tilde\varphi^{c,\eps}$ extended to $ B_{d_y}\left(\sqrt{2d_y}\|\mathcal{Y}\|\right)$, and again observe that $\tilde\varphi^{c,\eps}\big|_\cY=\varphi^{c,\eps}$. Following the proof of \cref{lemma:regularity_EOT}, one readily verifies that for any $s\geq 2$, we have
\begin{equation}
    \begin{aligned}
    \label{derivative_est_of_extension}
\|\tilde{\varphi}\|_{\infty,B_{d_x}\left(\sqrt{2d_x}\|\mathcal{X}\|\right)} & \leq \frac{5}{2}\|c\|_{\infty, B_{d_x}\left(\sqrt{2 d_x}\|\mathcal{X}\|\right)\times  \mathcal{Y}} \\
\left\|D^\alpha \tilde{\varphi}\right\|_{\infty,B_{d_x}\left(\sqrt{2d_x}\|\mathcal{X}\|\right)} & \leq C_s C_{c,\|\mathcal{X}\|,\|\mathcal{Y}\|}^s\left(1+\eps^{1-s}\right), \forall \alpha \in \mathbb{N}_0^{d_x},\ 1 \leq|\alpha| \leq s.
\end{aligned}
\end{equation}

Recall that $s_{\mathsf{KB}}=\left\lfloor d_x / 2\right\rfloor+3$ and set
\begin{equation}
    \label{eq:constant_C_dx_dy}
    C_{c,d_x,\mathcal{X},\mathcal{Y}}\coloneqq \left(\frac{5}{2}\|c\|_{\infty,B_{d_x}\left(\sqrt{2 d_x}\|\mathcal{X}\|\right)\times  \mathcal{Y} }\right)\vee\left(C_{c,\|\mathcal{X}\|,\|\mathcal{Y}\|}^{s_{\mathsf{KB}}}\right).
\end{equation}
By \eqref{derivative_est_of_extension}, we now have
\begin{equation}
\label{eq:eot_constant_b}
  \max _{\alpha:\,|\alpha| \leq s_{\mathsf{KB}}}\left\|D^\alpha \tilde{\varphi}\right\|_{\infty, B_{d_x}\left(\sqrt{2d_x}\|\mathcal{X}\|\right)} \leq C_{s}C_{c,\mathcal{X},\mathcal{Y},d_x}(1+\eps^{1-s_{\mathsf{KB}}}) \coloneqq b,
\end{equation}
and so $\tilde{\varphi} \in \cC_b^{s_{\mathsf{KB}}}\big(B_{d_x}\left(\sqrt{2 d_x}\|\mathcal{X}\|\right)\big)$.

Noting that $\cX\subset B_{d_x}\left(\sqrt{2 d_x}\|\mathcal{X}\|\right)$, by \cref{prop:approx_restated}, there exists $f \in \cF_{k,d_x}\left(\bar c_{b,d_x}\right)$ such that 
\begin{equation}
\label{eq:EOTpotential_approx}
  \left\|\varphi-f\right\|_{\infty,\cX}\lesssim \bar{c}_{b,d_x}d_x^{\frac 12} k^{-\frac 12},  
\end{equation}
where $\bar{c}_{b,d_x} = b\, \bar{c}_{d_x}$ and $\bar{c}_{d_x}$ is defined as (see \cite[Equation (A.15)]{sreekumar2022neural})
\begin{equation}
\label{eq:constant_c_bar}
    \begin{aligned}
\bar{c}_{d_x} & \coloneqq \left(\kappa_{d_x}{d_x}^{\frac{3}{2}} \vee 1\right)\pi^{\frac{d_x}{2}}  \\
& \times \Gamma\left(\frac{d_x}{2}+1\right)^{-1}(\operatorname{rad}(\cX)+1)^{d_x} 2^{s_{\mathsf{KB}}} d_x\left(\frac{1-{d_x}^{\frac{s_{\mathsf{KB}}}{2}}}{1-\sqrt{d_x}}+{d_x}^{\frac{s_{\mathsf{KB}}}{2}}\right) \max _{\|\alpha\|_1 \leq s_{\mathsf{KB}}}\left\|D^\alpha \Psi\right\|_{\infty, B_d(0.5)},
\end{aligned}
\end{equation}
with $\kappa_{d_x}^2\coloneqq \left(d_x+{d_x}^{(s_{\mathsf{KB}}-1)}\right) \int_{\mathbb{R}^d}\left(1+\|\omega\|^{2\left(s_{\mathsf{KB}}-2\right)}\right)^{-1} d \omega$, $\operatorname{rad}(\cX)=0.5 \sup _{x, x^{\prime} \in \cX}\left\|x-x^{\prime}\right\|$ and $\Psi(x) \propto \exp \left(-\frac{1}{0.5-\|x\|^2}\right) \mathbbm{1}_{\{\|x\|<0.5\}}$ as the canonical mollifier normalized to have unit mass. 

\medskip
Our last step is to lift the sup-norm neural approximation bound on the semi-dual potential from \eqref{eq:EOTpotential_approx} to a bound on the approximation error of the corresponding EOT cost. The following lemma is proven in \cref{appen:proof_neural_approx_reduction}.
\begin{lemma}[Neural approximation error reduction]
\label{lemma:neural_approx_reduction}
Fix $(\mu,\nu) \in \mathcal{P}(\cX) \times \mathcal{P}(\cY)$ and let $\varphi$ be the semi-dual EOT potential for $\OT^\eps_c(\mu,\nu)$ from \cref{lemma:regularity_EOT}. For any $f \in \cF_{k,d_x}\left(a\right)$, we have
\[
    \left|\OT^{\eps}_c(\mu,\nu)-\OT^{\eps}_{c,k,a}(\mu,\nu)\right| \leq  2\left\|\varphi-f\right\|_{\infty, \cX}.
\]
\end{lemma}

Note that by our assumption, $\mathcal{X} \subseteq[-1,1]^{d_x}$ and $\mathcal{Y} \subseteq[-1,1]^{d_y}$, the constant $\bar{c}_{b, d_x}$ only depends on $c,d_x,d_y$. Hence, setting $a=\bar{c}_{b, d_x}$ and combining \cref{lemma:neural_approx_reduction} 
 with \eqref{eq:EOTpotential_approx}, we obtain
\[
\left|\OT^{\eps}_{c,k,a}(\mu,\nu)-\OT^\eps_{c}(\mu,\nu)\right|\lesssim2\bar{c}_{b,d_x}d_x^{\frac 12} k^{-\frac 12}\lesssim_{c, , d_x, d_y}\left(1+\frac{1}{\eps^{\left\lfloor\frac{d_x}{2}\right\rfloor}+2}\right) k^{-\frac{1}{2}}.
\]

\subsection{Proof of Lemma \ref{lem:estimation_error_bound}}
\label{appen:proof_estimation_error}
\medskip

Set $\cF^{c,\eps}(a)\coloneqq \left\{f^{c,\eps}: f \in \cF_{k,d_x}(a)\right\}$, and first bound
\begin{align}\label{eq:eot_emp_error_decomposition}
    &\mathbb{E}\left[\left|\OT^\eps_{c ,k,a}(\mu, \nu)-\widehat\OT^\eps_{c , k,a}\left(X^n, Y^n\right)\right|\right]\nonumber\\
    & \leq\underbrace{n^{-\frac{1}{2}} \mathbb{E}\left[\sup_{f \in \cF_{k,d_x}(a)} \mspace{-3mu}n^{-\frac{1}{2}} \left| \sum_{i=1}^n \big(f(X_i)-\mathbb{E}_{\mu}[f]\big) \right|  \right]}_{(\mathrm{I})} +\mspace{-2mu} \underbrace{n^{-\frac{1}{2}}\mathbb{E}\left[\sup_{f \in \cF^{c,\eps}(a)}\mspace{-3mu}n^{-\frac{1}{2}}\left| \sum_{j=1}^n  \big(f(Y_j)-\mathbb{E}_{\nu}[f]\big) \right|  \right]}_{(\mathrm{II})}.
\end{align}
To control these expected suprema, we again require regularity of the involved function, as stated in the next lemma.

\begin{lemma}
\label{lemma:c_transform_regularity}
    Fix $c \in \cC^{\infty}$, the $(c,\eps)$-transform of NNs class $\cF_{k, d_x}(a)$ satisfies the following uniform smoothness properties:
\begin{equation}
\begin{aligned}
    \|f^{c,\eps}\|_{\infty,\cY} &\leq \|c\|_{\infty,\cX\times\cY}+3a(\|\cX\|+1)\\
    \left\|D^\alpha f^{c,\eps}\right\|_{\infty,\cY} &\leq C_sC_{c,\|\cX\|,\|\cY\|}\left(1+\eps^{1-s}\right)\quad, \forall \alpha \in \mathbb{N}_0^{d_x} \text { with } 1 \leq|\alpha| \leq s
\end{aligned}
\end{equation}
for any $s \geq 2$ and $f \in \cF_{k, d_x}(a)$. In particular, taking $c=c_{\bA}$, for any $\bA  \in \mathcal{D}_{M }$, there exists a constant $R_s$ that depends only on $s,d_x,d_y$, such that
\[
\max _{\alpha:|\alpha|_1 \leq s}\left\|D^\alpha f^{c_{\bA},\eps}\right\|_{\infty, \cY} \leq R_s(1+a)(1+\eps^{1-s})
\]
for any $s \geq 2$ and $f \in \cF_{k, d_x}(a)$.
\end{lemma}
The only difference between Lemmas \ref{lemma:regularity_EOT} and \ref{lemma:c_transform_regularity} is that here we consider the $(c,\eps)$-transform of NNs, rather than of dual EOT potentials. As our NNs are also compactly supported and bounded, the derivation of this result is all but identical to the proof of \cref{lemma:regularity_EOT}, and is therefore omitted to avoid repetition.

We proceed to bound Terms $(\mathrm{I})$ and $(\mathrm{II})$ from \eqref{eq:eot_emp_error_decomposition}. For the first, consider
\begin{align}&\mathbb{E}\left[\sup_{f \in \cF_{k,d_x}(a)} n^{-\frac{1}{2}} \left| \sum_{i=1}^n \left(f(X_i)-\mathbb{E}_{\mu}[f]\right) \right|  \right]\\
   & \stackrel{(a)}{\lesssim} \mathbb{E}\left[\int_0^{\infty} \sqrt{\log N\left(\delta, \cF_{k,d_x}(a),\|\cdot\|_{2, \mu_n}\right)} d \delta\right]\nonumber\\
    & \leq \int_0^{\infty} \sqrt{\sup _{\gamma \in \mathcal{P}(\cX)} \log N\left(\delta, \cF_{k,d_x}(a),\|\cdot\|_{2, \gamma}\right)} d \delta\nonumber\\
    & \stackrel{(b)}{=} \int_0^{6 a(\|\cX\|+1)} \sqrt{\sup _{\gamma \in \mathcal{P}(\cX)} \log N\left(\delta, \cF_{k,d_x}(a),\|\cdot\|_{2, \gamma}\right)} d \delta\nonumber\\
    & \lesssim  a(\|\cX\|+1)  \int_0^1 \sqrt{\sup _{\gamma \in \mathcal{P}(\cX)} \log N\left( 3 a(\|\cX\|+1)\delta, \cF_{k,d_x}(a),\|\cdot\|_{2, \gamma}\right)}d \delta\nonumber\\
    & \stackrel{(c)}{\lesssim}a(\|\cX\|+1)d_x^{\frac{3}{2}},\label{eq:term1_bound}
\end{align}
where:

(a) follows by \cite[Corollary 2.2.8]{van1996springer} since $n^{-\frac{1}{2}} \sum_{i=1}^n \sigma_i f\left(X_i\right)$, where $\{\sigma_i\}_{i=1}^n$ are i.i.d Rademacher random variables, is sub-Gaussian w.r.t. pseudo-metric $\|\cdot\|_{2, \mu_n}$ (by Hoeffding's inequality) ;

(b) is since $\sup _{x \in \mathcal{X}, f \in \cF_{k,d_x}(a)} |f(x)| \leq 3 a(\|\cX\|+1)$ and $N\left(\delta, \cF_{k,d_x}(a),\|\cdot\|_{2, \gamma}\right)=1$, whenever $\delta > 6 a(\|\cX\|+1)$;

(c) uses the bound $\int_0^1 \sqrt{\sup _{\gamma \in \mathcal{P}(\cX)} \log N\left(3 a(\|\cX\|+1) \delta, \cF_{k,d_x}(a),\|\cdot\|_{2, \gamma}\right)} d \delta \lesssim d_x^{\frac{3}{2}}$, which follows from step (A.33) in \cite{sreekumar2022neural}.

\medskip
For Term $(\mathrm{II})$, let $s=\left\lceil d_y  / 2\right\rceil+1$, and consider
\begin{align}
    \mathbb{E}\left[\sup _{f \in \cF^{c,\eps}(a)} n^{-\frac{1}{2}}\left|\sum_{i=1}^n\big(f\left(Y_i\right)-\mathbb{E}_\mu[f]\big)\right|\right] 
    &\stackrel{(a)}{\lesssim} \int_0^{6 a(\|\cX\|+1)} \sqrt{\sup _{\gamma \in \mathcal{P}(\cY)} \log N\left(\delta, \cF^{c,\eps}(a),\|\cdot\|_{2, \gamma}\right)} d \delta\nonumber\\
    &\lesssim \int_0^{6a(\|\cX\|+1)} \sqrt{\sup _{\gamma \in \mathcal{P}(\cY)} \log N_{[\,]}\left(2\delta, \cF^{c,\eps}(a),\|\cdot\|_{2, \gamma}\right)} d \delta\nonumber\\
    & \stackrel{(b)}{\lesssim} K_s \int_{0}^{6 a(\|\cX\|+1)} \left(\frac{R_s(1+a)(1+\eps^{1-s})}{2\delta}\right)^{\frac{d_y}{2s}}d\delta\nonumber\\
    & \lesssim_{c,d_x,d_y} a(1+a)(1+\eps^{1-s}),\label{eq:eot_term2_bound1}
\end{align}
where $K_s$ is a constant depends only on $s,d_x$ and $R_s$ is the constant from \cref{lemma:c_transform_regularity} (which depends on $c,s,d_x,d_y$ by $\mathcal{X} \subseteq[-1,1]^{d_x}$ and $\mathcal{Y} \subseteq[-1,1]^{d_y}$), (a) follows by a similar argument to that from the bound on Term $(\mathrm{I})$, along with equation \eqref{eq:lipshcitz_c_tranform}, which specifies the upper limit for entropy integral, while (b) follows by \cref{lemma:c_transform_regularity} and \cite[Corollary 2.7.2]{van1996springer}, which upper bounds the bracketing entropy number of smooth functions on a bounded convex support.

To arrive at the effective error bound from \cref{lem:estimation_error_bound}, we provide a second bound on Term $(\mathrm{II})$. This second bound yields a better dependence on dimension (namely, only the smaller dimension $d_x$ appears in the exponent) at the price of another $\sqrt{k}$ factor. Neither bound is uniformly superior over the other, and hence our final result will simply take the minimum of the two. By \eqref{eq:lipshcitz_c_tranform} from the proof of \cref{lemma:neural_approx_reduction}, we have 
\[
N\left(\delta, \mathcal{F}^{c,\eps}(a),\|\cdot\|_{2, \gamma}\right) \leq N\left(\delta, \cF_{k, d_x}(a),\|\cdot\|_{\infty,\cX}\right).
\]
Invoking Lemma 2 from \cite{sreekumar2022neural}, which upper bounds the metric entropy of ReLU NNs class on the RHS above, we further obtain
\begin{equation}
\label{eq:covering_bound_c_transform}
 \log N\left(\delta, \cF_{k, d_x}(a),\|\cdot\|_{\infty, \cX}\right) \leq \big((d_x+2) k+d_x+1\big)\log\left(1+20 a(\|\cX\|+1) \delta^{-1}\right),
\end{equation}
and proceed to bound Term $(\mathrm{II})$ as follows:
\begin{align}
     &\mathbb{E}\left[\sup_{f \in \mathcal{F}^{c,\eps}}n^{-\frac{1}{2}}\left| \sum_{j=1}^n  \big(f(Y_j)-\mathbb{E}_{\nu}[f]\big) \right|  \right]\\ &\lesssim a(\|\mathcal{X}\|+1)\int_0^1 \sqrt{\sup _{\gamma \in \mathcal{P}(\cX)} \log N\left(3a(\|\mathcal{X}\|+1) \delta, \mathcal{F}^{c,\eps}(a),\|\cdot\|_{2, \gamma}\right)} d \delta  \nonumber \\
     & \lesssim a(\|\mathcal{X}\|+1)d_x^{\frac{1}{2}}\sqrt{k}\int_0^1 \sqrt{\log(1+7\delta^{-1})} d \delta\nonumber\\
     & \lesssim a(\|\mathcal{X}\|+1)d_x^{\frac{1}{2}}\sqrt{k}.\label{eq:eot_term2_bound2}
\end{align}

Inserting \eqref{eq:term1_bound}, \eqref{eq:eot_term2_bound1}, and \eqref{eq:eot_term2_bound2} back into \eqref{eq:eot_emp_error_decomposition}, we obtain the desired bound on the empirical estimation error by setting $a=\bar{c}_{b,d_x}$, as was defined in approximation error analysis:
\begin{equation}
\begin{aligned}
    \label{eq:EOT_est_error}
    \mathbb{E}&\left[\left|\OT^\eps_{c ,k,a}(\mu, \nu)-\widehat\OT^\eps_{c , k,a}\left(X^n, Y^n\right)\right|\right]\\
    &\qquad\qquad\qquad\qquad\qquad\lesssim_{c, d_x, d_y} \min\left\{1+\frac{1}{\eps^{\left\lceil d_x+\frac{d_y}{2}\right\rceil+4}}\,,\left(1+\frac{1}{\eps^{\left\lfloor \frac{d_x}{2}\right\rfloor+2}}\right)\sqrt{k}\right\} n^{-\frac{1}{2}}.
\end{aligned}
\end{equation}

\subsection{Proof of Theorem \ref{thm:Bound_EOT_2}}
\label{appen:proof_Bound_EOT2}

This proof is similar to that of \cref{thm:Bound_EOT1}, up to minor modifications. For brevity, we only highlight the required changes. Note that for $k$ with $m_k \geq \bar{c}_{b,d_x}$, where the latter is defined in proof of \cref{thm:Bound_EGW1} (see \eqref{eq:eot_constant_b} and \eqref{eq:constant_c_bar}), we have $\cF_{k,d_x}(\bar{c}_{b,d_x})\subset\cF_{k,d_x}(m_k)$. Hence, by \cref{lemma:neural_approx_reduction} and \eqref{eq:EOTpotential_approx}, there exists a NN $f\in\cF_{k,d_x}(\bar{c}_{b,d_x})$, such that
\[
\left|\OT^{\eps}_{c,k,m_k}(\mu,\nu)-\OT^\eps_{c}(\mu,\nu)\right|\leq 2\left\|\varphi-f\right\|_{\infty,\cX}\lesssim_{c, d_x, d_y}\left(1+\frac{1}{\eps^{\left\lfloor\frac{d_x}{2}\right\rfloor}+2}\right) k^{-\frac{1}{2}}.
\]

Next, for estimation error, by setting $a=m_k\coloneqq\log k\vee 1$ in  \eqref{eq:term1_bound}, \eqref{eq:eot_term2_bound1}, and \eqref{eq:eot_term2_bound2} (instead of $a=\bar{c}_{b,d_x}$ as in the proof of \cref{lem:estimation_error_bound}), we arrive at
\begin{align*}
&\mathbb{E}\left[\left|\OT^\eps_{c ,k,m_k}(\mu, \nu)-\widehat\OT^\eps_{c , k,m_k}\left(X^n, Y^n\right)\right|\right] \\
&\qquad\qquad\qquad\qquad\qquad\qquad\lesssim_{c,  d_x, d_y}\min\left\{\left(1+\frac{1}{\eps^{\left[\frac{d_y}{2}\right\rceil}}\right)(\log k)^2\,, \sqrt{k}\log k\right\} n^{-\frac{1}{2}}.
\end{align*}
Combining the two errors yields:
\begin{align*}
   &\sup _{(\mu, \nu) \in \mathcal{P}(\cX) \times \mathcal{P}(\cY)} \mspace{-5mu}\mathbb{E}\left[\left|\widehat\OT^\eps_{c , k,m_k}\left(X^n, Y^n\right)-\mspace{-2mu}\OT_c^{\eps}(\mu, \nu)\right|\right]\\
   &\qquad\qquad\lesssim_{c,  d_x, d_y} \left(1+\frac{1}{\eps^{\left\lfloor\frac{d_x}{2}\right\rfloor+2}}\right) k^{-\frac{1}{2}}+\min\left\{\left(1+\frac{1}{\eps^{\left[\frac{d_y}{2}\right\rceil}}\right)(\log k)^2\,, \sqrt{k}\log k\right\} n^{-\frac{1}{2}}.
\end{align*}
\qed

\subsection{Proof of Theorem \ref{thm:Bound_NeuralCoupling_EOT}}
\label{appen:proof_NeuralCoupling_EOT}
Define $\Gamma(f): = \int_{\cX} fd\mu+\int_{\cY} f^{c, \eps}d \nu$, and let $\varphi_\star$ be optimal potential of $\OT_c^{\eps}(\mu,\nu)$, solving semi-dual formulation. Denote the corresponding optimal coupling by $\gamma^\eps_\star$. We first show that for any continuous $f:\cX\to\RR$, the following holds:
\begin{equation}
\label{eq:stability_coupling}
    \Gamma\left(\varphi_\star
    \right)-\Gamma\left(f\right) =\eps\KL\big(\gamma^\eps_\star\big\|\gamma^\eps_{f}\big),
\end{equation}
where
\[
d \gamma^\eps_{f}(x, y) = \frac{\exp \left(\frac{f(x)-c(x, y)}{\eps}\right)}{\int_{\cX} \exp \left(\frac{f-c(\cdot, y)}{\eps}\right) d\mu} d\mu\otimes\nu(x,y).
\]
The derivation is inspired by the proof of \cite[Theorem 2]{mokrov2023energy}, with several technical modifications. 
Since $\mu \in \mathcal{P}_{\mathsf{ac}}(\cX)$ with Lebesgue density $\frac{d\mu}{dx}$, define its energy function $E_{\mu}: \cX \rightarrow \mathbb{R}$ by $\frac{d \mu(x)}{d x} \propto \exp \left(-E_{\mu}(x)\right)$. Also define conditional distribution $d\gamma^\eps_{f}(\cdot| y) \coloneqq  \frac{d\pi_{f}(\cdot,y)}{d\nu(y)}$, and set $\tilde{f}\coloneqq f-\eps E_{\mu}(x)$. We have
\begin{align*}
    \frac{d\gamma^\eps_{f}(x| y) }{d x}& =\frac{d\gamma^\eps_{f}(x,y)}{d\nu(y)dx}\mspace{-.5mu}=\mspace{-.5mu} \frac{\exp \left(\frac{f(x)-c(x, y)}{\eps}\right)\frac{d\mu(x)}{dx}}{\int_{\cX} \exp \left(\frac{f(x')-c(x', y)}{\eps}\right)\frac{d\mu}{dx'}(x')dx'}\mspace{-.5mu}=\mspace{-.5mu} \frac{\exp \left(\frac{\tilde{f}(x)-c(x, y)}{\eps}\right)}{\int_{\cX} \exp \left(\frac{\tilde{f}(x')-c(x', y)}{\eps}\right)dx'}.
\end{align*}  
Define the shorthands $F_{f}(y)\coloneqq\int_{\cX} \exp \left(\frac{f(x)-c(x, y)}{\eps}\right)dx$ and $Z \coloneqq \int_\cX \exp \big(-E_{\mu}(x)\big)dx$, and note that the $(c,\eps)$-transform of $f$ can be expressed as
\begin{align*}
   f^{c, \eps}(y)&=-\eps \log \left(\int_{\cX} \exp \left(\frac{f(x)-c(x, y)}{\eps}\right) \frac{d \mu(x)}{dx}dx\right) \\
   &=-\eps \log \left(\int_{\cX} \exp \left(\frac{\tilde{f}(x)-c(x, y)}{\eps}\right)dx\right) + \eps\log(Z)\\
   &= -\eps \log \big(F_{\tilde{f}}(y)\big)+ \eps\log(Z).
\end{align*}
We are now ready to prove \eqref{eq:stability_coupling}. For $\rho\in\cP_{\mathsf{ac}}(\cX)$ with Lebesgue density $\frac{d\rho}{dx}$, denote the differential entropy of $\rho$ by $\sH(\rho)\coloneqq -\int_\cX \log\left(\frac{d\rho}{dx}\right)d\rho$. Consider:
\begin{align*}
    &\Gamma\left(\varphi_\star\right)-\Gamma\left(f\right)\\
    &=\int_{\cX \times \cY} c d \gamma^\eps_\star-\eps \int_{\cY} \sH\big(\gamma^\eps_\star(\cdot| y)\big) d \nu(y)+\eps \sH(\mu)-\int_{\cX} f d \mu-\int_{\cY} f^{c,\eps}d \nu\\
    &=\int_{\cX \times \cY}\big(c(x, y)-\tilde{f}(x)\big) d \gamma^\eps_\star(x, y)-\eps \int_{\cY} \sH\big(\gamma^\eps_\star(\cdot| y)\big) d \nu(y)+\eps \int_{\cY}\log\big(F_{\tilde f}\big)d\nu\\
     &=-\mspace{-2mu}\eps \mspace{-5mu}\int_{\cX \times \cY} \mspace{-22mu}\frac{\tilde{f}(x)-c(x, y)}{\eps} d \gamma^\eps_\star(x, y)\mspace{-2mu}+\mspace{-2mu}\eps \mspace{-5mu}\int_{\cX \times \cY} \mspace{-22mu}\log \big(F_{\tilde{f}}(y)\big) d\gamma^\eps_\star(x, y)-\eps \int_{\cY} \sH\big(\gamma^\eps_\star(\cdot| y)\big) d\nu(y)\\
     &=-\eps\int_{\cX \times \cY} \log \left(\frac{1}{F_{\tilde{f}}(y)} \exp \left(\frac{\tilde{f}(x)-c(x, y)}{\eps}\right)\right) d\gamma^\eps_\star(x, y)-\eps \int_{\cY} \sH\big(\gamma^\eps_\star(\cdot| y)\big) d\nu(y)\\
     &=-\eps\mspace{-4mu} \int_{\cX \times \cY} \mspace{-22mu}\log \left(\frac{d\gamma^\eps_{f}(x\mid y) }{dx}\right) d \gamma^\eps_\star(x, y)\mspace{-2mu}-\mspace{-2mu}\eps \mspace{-5mu}\int_{\cY} \sH\big(\gamma^\eps_\star(\cdot| y)\big) d\nu(y)
     \end{align*}
     \begin{align*}
     &=-\eps \int_{\cY}\int_{ \cX} \log \left(\frac{d\gamma^\eps_{f}(x| y) }{d x}\right) d \gamma^\eps_\star(x| y)d\nu(y)+\eps \int_{\cY}\int_{\cX} \log\left(\frac{d \gamma^\eps(x| y)}{d x}\right)d \gamma^\eps_\star(x| y) d\nu(y)\\
     &=\eps \int_{\cY}\int_{\cX}\log \left(\frac{d \gamma^\eps_\star(x | y)}{d\gamma^\eps_{f}(x|y)}\right)d\gamma^\eps_\star(x| y)d\nu(y)\\
     &=\eps\int_{\cX \times \cY} \log \left(\frac{d\gamma^\eps_\star(x, y)}{d\gamma^\eps_{f}(x, y)}\right) d\gamma^\eps_\star(x, y)\\
     &=\eps\KL\left(\gamma^\eps_\star\middle\|\gamma^\eps_{f}\right).
\end{align*}
Recalling that $\hat f_\star$ is a NN that optimizes the NE $\widehat\OT^\eps_{c,k,a}(X^n,Y^n)$ from \eqref{eq:NE_EOTc}, and plugging it into \eqref{eq:stability_coupling} yields $\KL\big(\pi^\eps_\star\big\|\pi^\eps_{\hat{f}_\star}\big)=\eps^{-1}\big[\Gamma_{} \left(\varphi_\star\right)-\Gamma(\hat{f}_\star)\big]$. Thus, to prove the KL divergence bound from Theorem \ref{thm:Bound_NeuralCoupling_EOT}, it suffices to control the gap between the $\Gamma$ functionals on the RHS above.

\medskip
Write $f_{\star}$ for a NN that maximizes the population-level neural EOT cost $\OT_{c,k,a}^\eps(\mu,\nu)$ (see \eqref{eq:neural_pop_EOT}). Define $\widehat{\Gamma}(f): = \frac 1n\sum_{i=1}^n\big[f(X_i)+f^{c, \eps}(Y_i)\big]$ for the optimization objective in the problem $\widehat\OT_{c,k,a}^\eps(X^n,Y^n)$ (see \eqref{eq:NE_EOTc}), and note that $\hat{f}_{\star}$ is a maximizer of $\widehat{\Gamma}$. We now have
\begin{equation}
\begin{aligned}
&\eps\KL\left(\gamma_\star^\eps\middle\|\gamma^\eps_{\hat{f}_{\star}}\right)\\
& = \Gamma\left(\varphi_{\star}\right)-\Gamma \left(f_{\star}\right)+\Gamma\left(f_\star\right)-\widehat{\Gamma} \left(\hat{f}_\star\right)+\widehat{\Gamma} \left(\hat{f}_\star\right)-\Gamma\left(\hat{f}_\star\right)\\
&=\underbrace{\OT_c^{\eps}(\mu, \nu)\mspace{-2mu}-\mspace{-2mu}\OT^{\eps}_{c,k,a}(\mu, \nu)}_{(\mathrm{I})} + \underbrace{\OT^{\eps}_{c, k,a}(\mu, \nu)\mspace{-2mu}-\mspace{-2mu}\widehat\OT^{\eps}_{c,k,a}\left(X^n, Y^n\right)}_{(\mathrm{II})}+\underbrace{\widehat{\Gamma} \left(\hat{f}_\star\right)\mspace{-2mu}-\mspace{-2mu}\Gamma \left(\hat{f}_\star\right)}_{(\mathrm{III})}.
\end{aligned}
\end{equation}
Setting $a =\bar{c}_{b,d_x}$ as in the proof of Theorem \ref{thm:Bound_EOT1} (see \eqref{eq:eot_constant_b} and \eqref{eq:constant_c_bar}) and taking an expectation (over the data) on both sides, Terms $(\mathrm{I})$ and $(\mathrm{II})$ are controlled, respectively, by the approximation error and empirical estimation error from \cref{lem:approximation_error_bound} and \cref{lem:estimation_error_bound}. For Term $(\mathrm{III})$, consider
\begin{align*}
    &\mathbb{E}\left[\widehat{\Gamma} \left(\hat{f}_\star\right)-\Gamma \left(\hat{f}_\star\right)\right]\\
    &\leq \mathbb{E}\left[\sup_{f \in \cF_{k, d_x}(\bar{c}_{b,d_x})} \left|\Gamma(f)-\widehat{\Gamma}(f)\right|\right]\\
    &\leq n^{-\frac{1}{2}}\mathbb{E}\left[\sup _{f \in \cF_{k, d_x}(\bar{c}_{b,d_x})}\mspace{-10mu} n^{-\frac{1}{2}}\left|\sum_{i=1}^n\big(f\left(X_i\right)-\mathbb{E}_\mu[f]\big)\right|\right]\mspace{-2mu}\mspace{-2mu}+\mspace{-2mu}n^{-\frac{1}{2}}\mathbb{E}\left[\sup _{f \in \mathcal{F}^{c, \eps}} \mspace{-3mu}n^{-\frac{1}{2}}\left|\sum_{j=1}^n\big(f\left(Y_j\right)-\mathbb{E}_\nu[f]\big)\right|\right].
\end{align*}
where the last step follows similarly to \eqref{eq:eot_emp_error_decomposition}. Notably, the RHS above is also bounded  by the estimation error bound from \eqref{eq:EGW_est_error}. Combining the above completes the proof.
\qed

\subsection{Proof of Theorem \ref{thm:Bound_EGW1}}
\label{appen:proof_Bound_EGW1}

For $(\mu, \nu) \in \mathcal{P}(\cX) \times \mathcal{P}(\cY)$, define the population-level neural EGW cost as
\[\mathsf{GW}^{\eps}_{k,a}(\mu, \nu)\coloneqq \sC(\mu,\nu)+\inf_{\bA \in\cD_{M }}\Big\{32\|\bA \|_\F^2+\OT^{\eps}_{\bA ,k,a}(\mu,\nu)\Big\},
\]
where we have used the shorthand $\OT^{\eps}_{\bA,k,a}(\mu, \nu)=\OT^{\eps}_{c_\bA,k,a}(\mu, \nu)$, which is the population-level neural EOT cost as defined in \eqref{eq:neural_pop_EOT} with cost function $c_\bA$. We decompose the neural estimation error into the approximation and empirical estimation errors:
\begin{align}
&\mathbb{E}\left[\left|\widehat{\mathsf{GW}}_{k,a}^{\mspace{1mu}\eps}\left(X^n, Y^n\right)-\mathsf{GW}^{\eps}(\mu, \nu)\right|\right]\\
&\leq\left|\mathsf{GW}_{k,a}^{\mspace{1mu}\eps}(\mu, \nu)-\mathsf{GW}^{\eps}(\mu, \nu)\right|+\mathbb{E}\left[\left|\mathsf{GW}^{\eps}_{k,a}(\mu, \nu)-\widehat{\mathsf{GW}}^{\eps}_{k,a}\left(X^n, Y^n\right)\right|\right] \nonumber\\
& \leq \underbrace{\sup _{\bA  \in \mathcal{D}_{M }}\left|\OT^{\eps}_{\bA ,k,a}(\mu, \nu)-\OT^{\eps}_{\bA}(\mu, \nu)\right|}_{\text {Approximation error}}+\underbrace{\mathbb{E}\left[\left|\sC(\mu, \nu)-\sC\left(X^n, Y^n\right)\right|\right]}_{\text {Estimation error 1}}\nonumber\\
&\qquad\qquad\qquad\qquad\qquad\qquad\qquad\qquad +\underbrace{\mathbb{E}\left[\sup _{\bA  \in \mathcal{D}_{M }} \left|\OT^{\eps}_{\bA , k,a}(\mu, \nu)-\widehat\OT^{\eps}_{\bA , k,a}\left(X^n, Y^n\right)\right|\right]}_{\text {Estimation error 2}},\label{eq:EGW_NE_error_decomposition}
\end{align}
and analyze each term separately. 

\medskip

\noindent\underline{Approximation error.} With specific expression of cost function $c_\bA$, we can have explicit bound on derivatives of potentials. From \eqref{eq:regularity_EGW_potential} in the proof of \cref{lemma:regularity_EOT}, we have that there exist semi-dual EOT potentials $(\varphi_{\bA},\varphi^{c,\eps}_{\bA})$ for $\OT_{\bA}^{\mspace{1mu}\eps}(\mu, \nu)$, such that
\begin{equation}
\begin{aligned}
\label{eq:regularity_EGW}
    \|\varphi_{\bA}\|_{\infty,\cX} &\leq 6 d_x d_y+40 (d_xd_y)^{3/2}\\
    \left\|D^\alpha \varphi_{\bA}\right\|_{\infty,\cX} &\leq C_{s}\left(1+\eps^{1-s}\right)\left(1+d_y\right)^{s}\left(1+\sqrt{d_x} d_y+\sqrt{d_x}\right)^{s},\quad \forall \alpha \in \mathbb{N}_0^{d_x} \text { with } 1 \leq|\alpha| \leq s,
\end{aligned}
\end{equation}
for any $s \geq 2$ and some constant $C_s$ that depends only on $s$. Note that $\|\calX\|\vee\|\calY\|\leq1$, we now have 
\begin{equation}
\label{eq:egw_constant_b}
  \max _{\alpha:\,|\alpha| \leq s_{\mathsf{KB}}}\left\|D^\alpha \tilde{\varphi}_{\bA}\right\|_{\infty, B_{d_x}\left(\sqrt{2d_x}\right)} \leq C_{s_{\mathsf{KB}}}C_{d_x,d_y}(1+\eps^{1-s_{\mathsf{KB}}}) \coloneqq b,
\end{equation}
where
\begin{equation}
    \label{eq:constant_C_dx_dy_egw}
    C_{d_x,d_y}\coloneqq \left\{\left(1+ d_y\right)^{s_{\mathsf{KB}}}\left(1+ \sqrt{d_x}d_y+\sqrt{ d_x}\right)^{s_{\mathsf{KB}}}\right\}\vee\left\{10 d_xd_y+(24+16\sqrt{2}) (d_x d_y)^{3/2} \right\}.
\end{equation}
Let $\bar{c}_{b, d_x}=b \bar{c}_{d_x}$ and set $a=\bar{c}_{b, d_x}$, where $\bar{c}_{d_x}$ is given in \eqref{eq:constant_c_bar}. Following \cref{lem:approximation_error_bound}, we obtain
\[
\left|\OT^{\eps}_{\mathbf{A},k,a}(\mu,\nu)-\OT^\eps_{\mathbf{A}}(\mu,\nu)\right|\lesssim 2\bar{c}_{b,d_x}d_x^{\frac{1}{2}}k^{-\frac{1}{2}}\lesssim_{d_x, d_y}\left(1+\frac{1}{\eps^{\left\lfloor \frac{d_x}{2}\right\rfloor+2}}\right) k^{-\frac{1}{2}}.
\]
As the RHS above is independent of $\mathbf{A}\in\cD_{M }$, we conclude that
\begin{equation}
    \label{eq:EGW_approx_error}
    \sup_{\mathbf{A}\in\cD_{M }} \left|\OT^{\eps}_{\mathbf{A},k,a}(\mu,\nu)-\OT^\eps_{\mathbf{A}}(\mu,\nu)\right|\lesssim_{d_x, d_y}\left(1+\frac{1}{\eps^{\left\lfloor \frac{d_x}{2}\right\rfloor+2}}\right) k^{-\frac{1}{2}}.
\end{equation}

\medskip
\noindent\underline{Estimation error 1.}
The estimation rate for $\mathsf{C}(\mu,\nu)$ was derived as part of the proof of \cite[Theorem 2]{zhang2024gromov} under a 4-sub-Weibull assumption on the population distributions.\footnote{A probability distribution $\mu\in\cP(\RR^d)$ is called $\beta$-sub-Weibull with parameter $\sigma^2$ for $\sigma\geq0$ if $\int \exp\left(\|x\|^\beta/2\sigma^2\right) d \mu(x)\leq 2$.} Since our $\mu,\nu$ are compactly supported, they are also 4-sub-Weilbull with parameter $\sigma^2=\frac{(d_x\vee d_y)^2}{2 \ln 2}$, and we consequently obtain
\begin{equation}
   \mathbb{E}\left[\left|\sC(\mu, \nu)-\sC\left(X^n, Y^n\right)\right|\right] \lesssim \frac{1+\sigma^4}{\sqrt{n}}.\label{eq:Sample_Complexity_C}
\end{equation}

\medskip

\noindent\underline{Estimation error 2}.  Set $\cF_{\mathcal{D}_{M }}^{c,\eps}(a)\coloneqq \bigcup_{\bA  \in \mathcal{D}_{M }}\left\{f^{c_{\bA},\eps}: f \in \cF_{k,d_x}(a)\right\}$, and first bound
\begin{align}\label{eq:emp_error_decomposition}
    &\mathbb{E}\left[\sup _{\bA  \in \mathcal{D}_{M }} \left|\OT^\eps_{\bA ,k,a}(\mu, \nu)-\widehat\OT^\eps_{\bA , k,a}\left(X^n, Y^n\right)\right|\right]\nonumber\\
    & \leq\underbrace{n^{-\frac{1}{2}} \mathbb{E}\left[\sup_{f \in \cF_{k,d_x}(a)} \mspace{-7mu}n^{-\frac{1}{2}} \left| \sum_{i=1}^n \big(f(X_i)-\mathbb{E}_{\mu}[f]\big) \right|  \right]}_{(\mathrm{I})} +\mspace{-2mu} \underbrace{n^{-\frac{1}{2}}\mathbb{E}\left[\sup_{f \in \cF_{\mathcal{D}_{M }}^{c,\eps}(a)}\mspace{-7mu}n^{-\frac{1}{2}}\left| \sum_{j=1}^n  \big(f(Y_j)-\mathbb{E}_{\nu}[f]\big) \right|  \right]}_{(\mathrm{II})}.
\end{align}
The bound for Term $(\mathrm{I})$ is given in \eqref{eq:term1_bound}. For Term  $(\mathrm{II})$, let $s=\left\lceil d_y / 2\right\rceil+1$, \eqref{eq:eot_term2_bound1} and \cref{lemma:c_transform_regularity} gives:
\begin{align}
    \mathbb{E}\left[\sup _{f \in \cF_{\mathcal{D}_{M }}^{c,\eps}(a)} n^{-\frac{1}{2}}\left|\sum_{i=1}^n\big(f\left(Y_i\right)-\mathbb{E}_\mu[f]\big)\right|\right] 
    & \lesssim K_s \int_{0}^{12 a} \left(\frac{R_s(1+a)(1+\eps^{1-s})}{2\delta}\right)^{\frac{d_y}{2s}}d\delta\nonumber\\
    & \lesssim_{d_x,d_y} a(1+a)(1+\eps^{1-s}),\label{eq:term2_bound1}
\end{align}

To arrive at the effective error bound from \cref{thm:Bound_EGW1}, we provide a second bound on Term $(\mathrm{II})$. However, we can not directly use the same analysis as in the proof of \cref{lem:estimation_error_bound}, since the function class here is a union of matrices $\bA$. Thus, we need to refine our analysis. For any $f,g\in\mathcal{F}_{k,d_x}$ and $\bA,\bB\in\mathcal{D}_M$, we have
\[
\left|f^{c_{\bA}, \eps}(y)-g^{c_{\bB}, \eps}(y)\right| \leq \left|f^{c_{\bA}, \eps}(y)-f^{c_{\bB}, \eps}(y)\right| +\left|f^{c_{\bB}, \eps}(y)-g^{c_{\bB}, \eps}(y)\right| 
\]
For the second term, by similar arguments to \eqref{eq:lipshcitz_c_tranform}, we have
\[
\left|f^{c_{\bB}, \eps}(y)-g^{c_{\bB}, \eps}(y)\right| \leq \left\|f-g\right\|_{\infty, \mathcal{X}}, \quad \forall y \in \mathcal{Y}
\]
For first term, by similar argument in proof of \eqref{eq:lipshcitz_c_tranform}, we have
\[
\left|f^{c_{\bA}, \eps}(y)-f^{c_{\bB}, \eps}(y)\right|  \leq 32\sqrt{d_xd_y}\|\calX\|\|\calY\|\|\bA-\bB\|_2\leq 32\sqrt{d_xd_y}\|\bA-\bB\|_2
\]
we obtain
\[
N(\delta,\mathcal{F}_{\mathcal{D}_M}^{c, \eps},\|\cdot\|_{2,\gamma}) \leq N\left(\frac{\delta}{64\sqrt{d_xd_y}},{\mathcal D}_M,||\cdot\|_2\right) \cdot N\left(\frac{\delta}{2}, \mathcal{F}_{k, d_x}(a),\|\cdot\|_{\infty, \mathcal{X}}\right)
\]
The covering number for ${\mathcal D}_M$ satisfies:
\[
N(\delta,{\mathcal D}_M,||\cdot\|_2) \leq \left(1+\frac{M \sqrt{d_x d_y}}{\delta}\right)^{d_x d_y}
\]
and \eqref{eq:covering_bound_c_transform} gives bound for $N\left(\delta, \mathcal{F}_{k, d_x}(a),\|\cdot\|_{\infty, \mathcal{X}}\right)$. Hence, we have
\begin{align}
     &\mathbb{E}\left[\sup_{f \in \mathcal{F}_{\mathcal{D}_M}^{c,\eps}}n^{-\frac{1}{2}}\left| \sum_{j=1}^n  \big(f(Y_j)-\mathbb{E}_{\nu}[f]\big) \right|  \right] \\
     &\lesssim a\int_0^1 \sqrt{\sup _{\gamma \in \mathcal{P}(\cX)} \log N\left(6 a \delta, \mathcal{F}_{\mathcal{D}_M}^{c,\eps}(a),\|\cdot\|_{2, \gamma}\right)} d \delta  \nonumber \\
     & \lesssim a\int_0^1 \sqrt{d_xd_y\log\left(1+\frac{32(d_xd_y)^{3/2}}{3a\delta}\right)+d_xk\log\left(1+14\delta^{-1}\right)} d \delta \nonumber\\
     & \lesssim a(d_xd_y)^{\frac{3}{2}}\sqrt{k} \sqrt{\frac{32\left(d_x d_y\right)^{3 / 2}}{3 a}+14}\nonumber\\
     &\lesssim a(d_xd_y)^{\frac{3}{2}}\sqrt{k}\label{eq:EGW_second_estimation_bound}
\end{align}
where the last step is due to our choice $3a \geq 32\left(d_x d_y\right)^{3 / 2}$.

\medskip

Inserting \eqref{eq:term1_bound}, \eqref{eq:term2_bound1}, and \eqref{eq:EGW_second_estimation_bound} back into \eqref{eq:emp_error_decomposition}, we obtain the desired bound on the empirical estimation error by setting $a=\bar{c}_{b,d_x}$, as was defined in approximation error analysis:
\begin{equation}
\begin{aligned}
    \label{eq:EGW_est_error}
    \mathbb{E}&\left[\sup _{\bA  \in \mathcal{D}_{M }} \left|\OT^\eps_{\bA ,k,a}(\mu, \nu)-\widehat\OT^\eps_{\bA , k,a}\left(X^n, Y^n\right)\right|\right]\\
    &\qquad\qquad\qquad\qquad\qquad\qquad\lesssim_{d_x,d_y} \min\left\{1+\frac{1}{\eps^{\left\lceil d_x+\frac{d_y}{2}\right\rceil+4}}\,,\left(1+\frac{1}{\eps^{\left\lfloor \frac{d_x}{2}\right\rfloor+2}}\right)\sqrt{k}\right\} n^{-\frac{1}{2}}.
\end{aligned}
\end{equation}

The proof is concluded by plugging the approximation error bound from \eqref{eq:EGW_approx_error} and the two estimation error bounds from \eqref{eq:Sample_Complexity_C} and \eqref{eq:EGW_est_error} into \eqref{eq:EGW_NE_error_decomposition}, and supremizing over $(\mu, \nu) \in \mathcal{P}(\cX) \times \mathcal{P}(\cY)$, while noting that all the above bounds holds uniformly in the two distributions. 
\qed

\section{Concluding Remarks and Outlook}

This work proposed a neural estimation framework subject to non-asymptotic accuracy guarantees for the EOT and the quadratic EGW problems between distributions supported on Euclidean spaces. The EOT NE can handle any smooth cost function and employed the semi-dual form to calibrate the estimator to the lower-dimensional space. 
The EOT estimate served as a module in the EGW estimation algorithm which leveraged the variational representation from \cite{zhang2024gromov} to express EGW as an infimum of a class of EOT problems parametrized by their costs. Accelerated first-order methods were used to optimize over the cost parameters, with each gradient computation relying on neural estimation of the current EOT cost.

Our approach yielded estimates not only for the EOT and EGW costs but also for the optimal transport/alignment plan. Non-asymptotic formal guarantees on the quality of the cost and plan NE were provided, under the sole assumption of compactly supported population distributions, with no further regularity conditions imposed. Our bounds revealed optimal scaling laws for the NN and the dataset sizes that ensure parametric (and hence minimax-rate optimal) convergence. In terms of dependence on the regularization parameter $\eps$ and data dimensions $d_x,d_y$, our bounds exhibited the LCA principle \cite{hundrieser2022empirical}, as they adapt to the smaller of the two dimensions. The proposed EGW NE was tested via numerical experiments on synthetic and real-world data, demonstrating its accuracy, scalability, and fast convergence rates that match the derived theory. 

 Future research directions stemming from this work are abundant. First, one may consider the neural estimation framework in the mean-field regime, i.e., when the number of neurons tens to infinity. In this over-parameterized setting, the training dynamics can be described as a gradient flow on the space of probability measures, a viewpoint developed through the lenses of distributional dynamics and measure optimization \cite{mei2018mean, chizat2018global}. This perspective shifts the analytical focus from the approximation–estimation trade-off, which is at the core of this work, to the interplay between optimization and statistical estimation (as the approximation error vanishes in the over-parametrized limit). Crucially, global convergence of the mean-field dynamics can be established for objectives, such as ours, that are convex functionals of the parameter measure. Developing this program for neural OT estimators, including an account of discretization, stability, and how this theory pairs with the outer EGW update scheme from \cref{alg:EGW_NE}, is an exciting direction for future work.

On a more technical level, our analysis accounts only for compactly supported distributions. It seems possible to extend our results to distributions with unbounded supports using the technique from \cite{sreekumar2022neural}that considers a sequence of restrictions to balls of increasing radii, along with the regularity theory for dual potentials from \cite{zhang2024gromov}, which accounts for 4-sub-Weibull distributions. Unfortunately, as in \cite{sreekumar2022neural}, rate bounds obtained from this technique would be sub-optimal. Obtaining sharp rates for the unboundedly supported case would require new ideas and forms an interesting research direction. 

Lastly, while EGW serves as an important approximation of GW, neural estimation of the GW distance itself is a challenging and appealing research avenue. The EGW variational representation from \cite{zhang2024gromov} specializes to the unregularized GW case by setting $\eps=0$, yielding a minimax objective akin to \eqref{eq:alg_objective}, but with the classical OT (Kantorovich-Rubinstein) dual instead of the EOT semi-dual. One may attempt to directly approximate this objective by neural networks, but dual OT potential generally lack sufficient regularity to allow quantitative approximation bounds. Assuming smoothness of the population distributions, and employing estimators that adapt to this smoothness, e.g., based on kernel density estimators or wavelets \cite{deb2021rates,manole2021plugin}, may enable deriving optimal convergence rates in the so-called \emph{high-smoothness regime}.

\bibliographystyle{unsrt}
\bibliography{ref}

@inproceedings{Wang2024neural,
  title={Neural estimation of entropic optimal transport},
  author = {Wang, Tao and Goldfeld, Ziv},
  booktitle={IEEE International Symposium on Information Theory (ISIT-2024)},
  year={2024}
}

@article{demetci2020gromov,
  title={Gromov-Wasserstein optimal transport to align single-cell multi-omics data},
  author={Demetci, Pinar and Santorella, Rebecca and Sandstede, Bj{\"o}rn and Noble, William Stafford and Singh, Ritambhara},
  journal={BioRxiv},
  pages={2020--04},
  year={2020},
  publisher={Cold Spring Harbor Laboratory}
}

@article{groppe2023lower,
  title={Lower Complexity Adaptation for Empirical Entropic Optimal Transport},
  author={Groppe, Michel and Hundrieser, Shayan},
  journal={arXiv preprint arXiv:2306.13580},
  year={2023}
}

@article{tsur2023data,
  title={Data-Driven Optimization of Directed Information over Discrete Alphabets},
  author={Tsur, Dor and Aharoni, Ziv and Goldfeld, Ziv and Permuter, Haim},
  journal={Accepted to the IEEE Transactions on Information theory},
  month = {November},
  year={2023}
}

@article{tsur2023neural,
  title={Neural estimation and optimization of directed information over continuous spaces},
  author={Tsur, Dor and Aharoni, Ziv and Goldfeld, Ziv and Permuter, Haim},
  journal={IEEE Transactions on Information Theory},
  year={2023},
  publisher={IEEE}
}

@inproceedings{gulrajani2017improved,
  title={Improved training of {W}asserstein {GAN}s},
  author={I. Gulrajani and F. Ahmed and M. Arjovsky and V. Dumoulin and A. C. Courville},
  booktitle={Proceedings of the Annual Conference on Advances in Neural Information Processing Systems (NeurIPS-2017)},
  pages={5767--5777},
  year={2017},
  month={Dec.},
  address={Long Beach, CA, US}
}

@inproceedings{arjovsky2017wasserstein,
 title={Wasserstein generative adversarial networks},
 author={M. Arjovsky and S. Chintala and L. Bottou},
 booktitle={Proceedings of the 34th International Conference on Machine Learning},
 pages={214--223},
 year={2017},
 month={Jul.},
 address={Sydney, Australia}
}

@article{beinert2023assignment,
  title={On assignment problems related to Gromov--Wasserstein distances on the real line},
  author={Beinert, Robert and Heiss, Cosmas and Steidl, Gabriele},
  journal={SIAM Journal on Imaging Sciences},
  volume={16},
  number={2},
  pages={1028--1032},
  year={2023},
  publisher={SIAM}
}

@article{mokrov2023energy,
  title={Energy-guided Entropic Neural Optimal Transport},
  author={Mokrov, Petr and Korotin, Alexander and Burnaev, Evgeny},
  journal={arXiv preprint arXiv:2304.06094},
  year={2023}
}

@article{daniels2021score,
  title={Score-based generative neural networks for large-scale optimal transport},
  author={Daniels, Max and Maunu, Tyler and Hand, Paul},
  journal={Advances in neural information processing systems},
  volume={34},
  pages={12955--12965},
  year={2021}
}

@inproceedings{sreekumar2021non,
  title={Non-asymptotic Performance Guarantees for Neural Estimation of $f$-Divergences},
  author={S. Sreekumar, Z. Zhang and Z. Goldfeld},
  booktitle = {International Conference on Artificial Intelligence and Statistics (AISTATS-2021)},
  series    = {Proceedings of Machine Learning Research},
  volume    = {130},
  pages={3322--3330},
  year      = {2021},
  month = {April},
  address = {Virtual conference}
}

@inproceedings{poole2018variational,
  title={On variational lower bounds of mutual information},
  author={Poole, Ben and Ozair, Sherjil and van den Oord, A{\"a}ron and Alemi, Alexander A and Tucker, George},
  booktitle={NeurIPS Workshop on Bayesian Deep Learning},
  year={2018}
}

@inproceedings{song2019understanding,
  title={Understanding the Limitations of Variational Mutual Information Estimators},
  author={Song, Jiaming and Ermon, Stefano},
  booktitle={International Conference on Learning Representations},
  year={2019}
}

@article{chan2019neural,
  title={Neural entropic estimation: A faster path to mutual information estimation},
  author={Chan, Chung and Al-Bashabsheh, Ali and Huang, Hing Pang and Lim, Michael and Tam, Da Sun Handason and Zhao, Chao},
  journal={arXiv preprint arXiv:1905.12957},
  year={2019}
}

@inproceedings{mroueh2021improved,
  title={Improved mutual information estimation},
  author={Mroueh, Youssef and Melnyk, Igor and Dognin, Pierre and Ross, Jarret and Sercu, Tom},
  booktitle={Proceedings of the AAAI Conference on Artificial Intelligence},
  volume={35},
  number={10},
  pages={9009--9017},
  year={2021}
}

@article{yukich1995sup,
  title={Sup-norm approximation bounds for networks through probabilistic methods},
  author={Yukich, Joseph E and Stinchcombe, Maxwell B and White, Halbert},
  journal={IEEE Transactions on Information Theory},
  volume={41},
  number={4},
  pages={1021--1027},
  year={1995},
  publisher={IEEE}
}

@article{goldfeld2022k,
  title={$k$-sliced mutual information: {A} quantitative study of scalability with dimension},
  author={Goldfeld, Ziv and Greenewald, Kristjan and Nuradha, Theshani and Reeves, Galen},
  journal={Advances in Neural Information Processing Systems},
  volume = {35},
  pages = {15982--15995},
  year={2022}
}

@article{tsur2023max,
  title={Max-Sliced Mutual Information},
  author={Tsur, Dor and Goldfeld, Ziv and Greenewald, Kristjan},
  journal={arXiv preprint arXiv:2309.16200},
  year={2023}
}

@article{molavipour2021neural,
  title={Neural Estimator of Information for Time-Series Data with Dependency},
  author={Molavipour, Sina and Ghourchian, Hamid and Bassi, Germ{\'a}n and Skoglund, Mikael},
  journal={Entropy},
  volume={23},
  number={6},
  pages={641},
  year={2021},
  publisher={Multidisciplinary Digital Publishing Institute}
}

@article{sreekumar2022neural,
  title={Neural estimation of statistical divergences},
  author={Sreekumar, Sreejith and Goldfeld, Ziv},
  journal={Journal of Machine Learning Research},
  volume={23},
  number={126},
  pages={1--75},
  year={2022}
}

@article{blumberg2020mrec,
  title={MREC: a fast and versatile framework for aligning and matching point clouds with applications to single cell molecular data},
  author={Blumberg, Andrew J and Carriere, Mathieu and Mandell, Michael A and Rabadan, Raul and Villar, Soledad},
  journal={arXiv preprint arXiv:2001.01666},
  year={2020}
}

@article{alvarez2018gromov,
  title={Gromov-Wasserstein alignment of word embedding spaces},
  author={Alvarez-Melis, David and Jaakkola, Tommi S},
  journal={arXiv preprint arXiv:1809.00013},
  year={2018}
}

@article{koehl2023computing,
  title={Computing the Gromov-Wasserstein Distance between Two Surface Meshes Using Optimal Transport},
  author={Koehl, Patrice and Delarue, Marc and Orland, Henri},
  journal={Algorithms},
  volume={16},
  number={3},
  pages={131},
  year={2023},
  publisher={MDPI}
}

@inproceedings{memoli2009spectral,
  title={Spectral Gromov-Wasserstein distances for shape matching},
  author={M{\'e}moli, Facundo},
  booktitle={2009 IEEE 12th International Conference on Computer Vision Workshops, ICCV Workshops},
  pages={256--263},
  year={2009},
  organization={IEEE}
}

@inproceedings{xu2019gromov,
  title={Gromov-wasserstein learning for graph matching and node embedding},
  author={Xu, Hongteng and Luo, Dixin and Zha, Hongyuan and Duke, Lawrence Carin},
  booktitle={International conference on machine learning},
  pages={6932--6941},
  year={2019},
  organization={PMLR}
}

@inproceedings{bunne2019learning,
  title={Learning generative models across incomparable spaces},
  author={Bunne, Charlotte and Alvarez-Melis, David and Krause, Andreas and Jegelka, Stefanie},
  booktitle={International conference on machine learning},
  pages={851--861},
  year={2019},
  organization={PMLR}
}

@article{mena2019statistical,
  title={Statistical bounds for entropic optimal transport: sample complexity and the central limit theorem},
  author={Mena, Gonzalo and Niles-Weed, Jonathan},
  journal={Advances in Neural Information Processing Systems},
  volume={32},
  year={2019}
}

@inproceedings{genevay2019sample,
  title={Sample complexity of sinkhorn divergences},
  author={Genevay, Aude and Chizat, L{\'e}naic and Bach, Francis and Cuturi, Marco and Peyr{\'e}, Gabriel},
  booktitle={The 22nd international conference on artificial intelligence and statistics},
  pages={1574--1583},
  year={2019},
  organization={PMLR}
}

@article{schmidt2020nonparametric,
  title={Nonparametric regression using deep neural networks with ReLU activation function},
  author={Schmidt-Hieber, Johannes},
  year={2020}
}

@article{bresler2020sharp,
  title={Sharp representation theorems for relu networks with precise dependence on depth},
  author={Bresler, Guy and Nagaraj, Dheeraj},
  journal={Advances in Neural Information Processing Systems},
  volume={33},
  pages={10697--10706},
  year={2020}
}

@article{shen2019deep,
  title={Deep network approximation characterized by number of neurons},
  author={Shen, Zuowei and Yang, Haizhao and Zhang, Shijun},
  journal={arXiv preprint arXiv:1906.05497},
  year={2019}
}

@article{deb2021rates,
  title={Rates of estimation of optimal transport maps using plug-in estimators via barycentric projections},
  author={Deb, Nabarun and Ghosal, Promit and Sen, Bodhisattva},
  journal={Advances in Neural Information Processing Systems},
  volume={34},
  pages={29736--29753},
  year={2021}
}

@article{manole2021plugin,
  title={Plugin estimation of smooth optimal transport maps},
  author={Manole, Tudor and Balakrishnan, Sivaraman and Niles-Weed, Jonathan and Wasserman, Larry},
  journal={arXiv preprint arXiv:2107.12364},
  year={2021}
}

@article{constantine1996multivariate,
  title={A multivariate Faa di Bruno formula with applications},
  author={Constantine, G and Savits, T},
  journal={Transactions of the American Mathematical Society},
  volume={348},
  number={2},
  pages={503--520},
  year={1996}
}

@article{memoli2011gromov,
  title={Gromov--Wasserstein distances and the metric approach to object matching},
  author={M{\'e}moli, Facundo},
  journal={Foundations of computational mathematics},
  volume={11},
  pages={417--487},
  year={2011},
  publisher={Springer}
}

@book{nesterov2003introductory,
  title={Introductory lectures on convex optimization: A basic course},
  author={Nesterov, Yurii},
  volume={87},
  year={2003},
  publisher={Springer Science \& Business Media}
}

@book{villani2008optimal,
    title={Optimal Transport: Old and New},
    author={C\'{e}dric Villani},
    year={2008},
    publisher={Springer}
}

@book{santambrogio2010,
author = {Santambrogio, Filippo},
title = {Optimal Transport for Applied Mathematicians},
publisher = {Birkh\"{a}user},
year = {2015}
}

@article{peyre2017computational,
  title={Computational optimal transport},
  author={Peyr{\'e}, Gabriel and Cuturi, Marco and others},
  journal={Center for Research in Economics and Statistics Working Papers},
  number={2017-86},
  year={2017}
}

@article{fournier2015rate,
  title={On the rate of convergence in Wasserstein distance of the empirical measure},
  author={Fournier, Nicolas and Guillin, Arnaud},
  journal={Probability theory and related fields},
  volume={162},
  number={3-4},
  pages={707--738},
  year={2015},
  publisher={Springer}
}

@article{xu2019scalable,
  title={Scalable Gromov-Wasserstein learning for graph partitioning and matching},
  author={Xu, Hongteng and Luo, Dixin and Carin, Lawrence},
  journal={Advances in neural information processing systems},
  volume={32},
  year={2019}
}

@article{tolstikhin2017wasserstein,
  title={Wasserstein auto-encoders},
  author={Tolstikhin, Ilya and Bousquet, Olivier and Gelly, Sylvain and Schoelkopf, Bernhard},
  journal={arXiv preprint arXiv:1711.01558},
  year={2017}
}

@article{courty2017joint,
  title={Joint distribution optimal transportation for domain adaptation},
  author={Courty, Nicolas and Flamary, R{\'e}mi and Habrard, Amaury and Rakotomamonjy, Alain},
  journal={Advances in neural information processing systems},
  volume={30},
  year={2017}
}

@inproceedings{chen2020graph,
  title={Graph optimal transport for cross-domain alignment},
  author={Chen, Liqun and Gan, Zhe and Cheng, Yu and Li, Linjie and Carin, Lawrence and Liu, Jingjing},
  booktitle={International Conference on Machine Learning},
  pages={1542--1553},
  year={2020},
  organization={PMLR}
}

@article{petric2019got,
  title={GOT: an optimal transport framework for graph comparison},
  author={Petric Maretic, Hermina and El Gheche, Mireille and Chierchia, Giovanni and Frossard, Pascal},
  journal={Advances in Neural Information Processing Systems},
  volume={32},
  year={2019}
}

@article{goldfeld2022limit,
  title={Limit theorems for entropic optimal transport maps and the Sinkhorn divergence},
  author={Goldfeld, Ziv and Kato, Kengo and Rioux, Gabriel and Sadhu, Ritwik},
  journal={arXiv preprint arXiv:2207.08683},
  year={2022}
}

@article{zhang2024gromov,
  title={Gromov--Wasserstein distances: Entropic regularization, duality and sample complexity},
  author={Zhang, Zhengxin and Goldfeld, Ziv and Mroueh, Youssef and Sriperumbudur, Bharath K},
  journal={The Annals of Statistics},
  volume={52},
  number={4},
  pages={1616--1645},
  year={2024},
  publisher={Institute of Mathematical Statistics}
}

@article{nguyen2010estimating,
  title={Estimating divergence functionals and the likelihood ratio by convex risk minimization},
  author={Nguyen, XuanLong and Wainwright, Martin J and Jordan, Michael I},
  journal={IEEE Transactions on Information Theory},
  volume={56},
  number={11},
  pages={5847--5861},
  year={2010},
  publisher={IEEE}
}

@inproceedings{le2022entropic,
  title={Entropic gromov-wasserstein between gaussian distributions},
  author={Le, Khang and Le, Dung Q and Nguyen, Huy and Do, Dat and Pham, Tung and Ho, Nhat},
  booktitle={International Conference on Machine Learning},
  pages={12164--12203},
  year={2022},
  organization={PMLR}
}

@inproceedings{seguy2018large,
  title={Large-Scale Optimal Transport and Mapping Estimation},
  author={Seguy, Vivien and Damodaran, Bharath Bhushan and Flamary, Remi and Courty, Nicolas and Rolet, Antoine and Blondel, Mathieu},
  booktitle={ICLR 2018-International Conference on Learning Representations},
  pages={1--15},
  year={2018}
}

@article{solomon2015convolutional,
  title={Convolutional {W}asserstein distances: {E}fficient optimal transportation on geometric domains},
  author={Solomon, Justin and De Goes, Fernando and Peyr{\'e}, Gabriel and Cuturi, Marco and Butscher, Adrian and Nguyen, Andy and Du, Tao and Guibas, Leonidas},
  journal={ACM Transactions on Graphics (ToG)},
  volume={34},
  number={4},
  pages={1--11},
  year={2015},
  publisher={ACM New York, NY, USA}
}

@article{li2020continuous,
  title={Continuous regularized {W}asserstein barycenters},
  author={Li, Lingxiao and Genevay, Aude and Yurochkin, Mikhail and Solomon, Justin M},
  journal={Advances in Neural Information Processing Systems},
  volume={33},
  pages={17755--17765},
  year={2020}
}

@article{zhang2024gradient,
  title={Gradient Flows and Riemannian Structure in the Gromov-Wasserstein Geometry},
  author={Zhang, Zhengxin and Goldfeld, Ziv and Greenewald, Kristjan and Mroueh, Youssef and Sriperumbudur, Bharath K},
  journal={arXiv preprint arXiv:2407.11800},
  year={2024}
}

@article{schiebinger2017reconstruction,
  title={Reconstruction of developmental landscapes by optimal-transport analysis of single-cell gene expression sheds light on cellular reprogramming},
  author={Schiebinger, Geoffrey and Shu, Jian and Tabaka, Marcin and Cleary, Brian and Subramanian, Vidya and Solomon, Aryeh and Liu, Siyan and Lin, Stacie and Berube, Peter and Lee, Lia and others},
  journal={BioRxiv},
  pages={191056},
  year={2017},
  publisher={Cold Spring Harbor Laboratory}
}

@article{bunne2024optimal,
  title={Optimal transport for single-cell and spatial omics},
  author={Bunne, Charlotte and Schiebinger, Geoffrey and Krause, Andreas and Regev, Aviv and Cuturi, Marco},
  journal={Nature Reviews Methods Primers},
  volume={4},
  number={1},
  pages={58},
  year={2024},
  publisher={Nature Publishing Group UK London}
}

@article{rioux2024limit,
  title={Limit Laws for {G}romov-{W}Wasserstein Alignment with Applications to Testing Graph Isomorphisms},
  author={Rioux, Gabriel and Goldfeld, Ziv and Kato, Kengo},
  journal={arXiv preprint arXiv:2410.18006},
  year={2024}
}

@article{repasky2023neural,
  title={Neural {Stein} Critics With Staged $L^2$-Regularization},
  author={Repasky, Matthew and Cheng, Xiuyuan and Xie, Yao},
  journal={IEEE Transactions on Information Theory},
  year={2023},
  publisher={IEEE}
}

@article{diaz2021lower,
  title={Lower bounds for the {MMSE} via neural network estimation and their applications to privacy},
  author={Diaz, Mario and Kairouz, Peter and Sankar, Lalitha},
  journal={arXiv preprint arXiv:2108.12851},
  year={2021}
}

@article{rioux2024entropic,
  title={Entropic {G}romov-{W}asserstein Distances: Stability and Algorithms},
  author={Rioux, Gabriel and Goldfeld, Ziv and Kato, Kengo},
  journal={Journal of Machine Learning Research},
  volume={25},
  number={363},
  pages={1--52},
  year={2024}
}

@article{barron1993universal,
  title={Universal approximation bounds for superpositions of a sigmoidal function},
  author={Barron, Andrew R},
  journal={IEEE Transactions on Information theory},
  volume={39},
  number={3},
  pages={930--945},
  year={1993},
  publisher={IEEE}
}

@inproceedings{barron1992neural,
  title={Neural net approximation},
  author={Barron, Andrew R},
  booktitle={Proc. 7th Yale workshop on adaptive and learning systems},
  volume={1},
  pages={69--72},
  year={1992}
}

@article{commander2005survey,
  title={A survey of the quadratic assignment problem, with applications},
  author={Commander, Clayton W},
  year={2005}
}

@article{solomon2016entropic,
  title={Entropic metric alignment for correspondence problems},
  author={Solomon, Justin and Peyr{\'e}, Gabriel and Kim, Vladimir G and Sra, Suvrit},
  journal={ACM Transactions on Graphics (ToG)},
  volume={35},
  number={4},
  pages={1--13},
  year={2016},
  publisher={ACM New York, NY, USA}
}

@inproceedings{peyre2016gromov,
  title={Gromov-wasserstein averaging of kernel and distance matrices},
  author={Peyr{\'e}, Gabriel and Cuturi, Marco and Solomon, Justin},
  booktitle={International conference on machine learning},
  pages={2664--2672},
  year={2016},
  organization={PMLR}
}

@article{cuturi2013sinkhorn,
  title={Sinkhorn distances: Lightspeed computation of optimal transport},
  author={Cuturi, Marco},
  journal={Advances in neural information processing systems},
  volume={26},
  year={2013}
}

@article{van1996springer,
  title={Springer series in statistics},
  author={van der Vaart, Aad W and Wellner, Jon A},
  journal={Weak convergence and empirical processesSpringer, New York},
  year={1996}
}

@article{vayer2019sliced,
  title={Sliced gromov-wasserstein},
  author={Vayer, Titouan and Flamary, R{\'e}mi and Tavenard, Romain and Chapel, Laetitia and Courty, Nicolas},
  journal={arXiv preprint arXiv:1905.10124},
  year={2019}
}

@article{belghazi2018mine,
  title={Mine: mutual information neural estimation},
  author={Belghazi, Mohamed Ishmael and Baratin, Aristide and Rajeswar, Sai and Ozair, Sherjil and Bengio, Yoshua and Courville, Aaron and Hjelm, R Devon},
  journal={arXiv preprint arXiv:1801.04062},
  year={2018}
}

@article{sejourne2021unbalanced,
  title={The unbalanced gromov wasserstein distance: Conic formulation and relaxation},
  author={S{\'e}journ{\'e}, Thibault and Vialard, Fran{\c{c}}ois-Xavier and Peyr{\'e}, Gabriel},
  journal={Advances in Neural Information Processing Systems},
  volume={34},
  pages={8766--8779},
  year={2021}
}

@inproceedings{scetbon2022linear,
  title={Linear-time gromov wasserstein distances using low rank couplings and costs},
  author={Scetbon, Meyer and Peyr{\'e}, Gabriel and Cuturi, Marco},
  booktitle={International Conference on Machine Learning},
  pages={19347--19365},
  year={2022},
  organization={PMLR}
}

@article{dumont2022existence,
  title={On the Existence of Monge Maps for the Gromov-Wasserstein Distance},
  author={Dumont, Theo and Lacombe, Th{\'e}o and Vialard, Fran{\c{c}}ois-Xavier},
  journal={arXiv preprint arXiv:2210.11945},
  year={2022}
}

@inproceedings{zhang2022cycle,
  title={Cycle consistent probability divergences across different spaces},
  author={Zhang, Zhengxin and Mroueh, Youssef and Goldfeld, Ziv and Sriperumbudur, Bharath},
  booktitle={International Conference on Artificial Intelligence and Statistics},
  pages={7257--7285},
  year={2022},
  organization={PMLR}
}

@article{barron1994approximation,
  title={Approximation and estimation bounds for artificial neural networks},
  author={Barron, Andrew R},
  journal={Machine learning},
  volume={14},
  pages={115--133},
  year={1994},
  publisher={Springer}
}

@article{bach2017breaking,
  title={Breaking the curse of dimensionality with convex neural networks},
  author={Bach, Francis},
  journal={The Journal of Machine Learning Research},
  volume={18},
  number={1},
  pages={629--681},
  year={2017},
  publisher={JMLR. org}
}

@article{suzuki2018adaptivity,
  title={Adaptivity of deep ReLU network for learning in Besov and mixed smooth Besov spaces: optimal rate and curse of dimensionality},
  author={Suzuki, Taiji},
  journal={arXiv preprint arXiv:1810.08033},
  year={2018}
}

@article{yang1999information,
  title={Information-theoretic determination of minimax rates of convergence},
  author={Yang, Yuhong and Barron, Andrew},
  journal={Annals of Statistics},
  pages={1564--1599},
  year={1999},
  publisher={JSTOR}
}

@article{uppal2019nonparametric,
  title={Nonparametric density estimation \& convergence rates for gans under besov ipm losses},
  author={Uppal, Ananya and Singh, Shashank and P{\'o}czos, Barnab{\'a}s},
  journal={Advances in neural information processing systems},
  volume={32},
  year={2019}
}

@article{hundrieser2022empirical,
  title={Empirical optimal transport between different measures adapts to lower complexity},
  author={Hundrieser, Shayan and Staudt, Thomas and Munk, Axel},
  journal={arXiv preprint arXiv:2202.10434},
  year={2022}
}

@article{kingma2014adam,
  title={Adam: A method for stochastic optimization},
  author={Kingma, Diederik P and Ba, Jimmy},
  journal={arXiv preprint arXiv:1412.6980},
  year={2014}
}

@article{birman1967piecewise,
  title={Piecewise-polynomial approximations of functions of the classes},
  author={Birman, M {\v{S}} and Solomjak, Mikhail Z},
  journal={Mathematics of the USSR-Sbornik},
  volume={2},
  number={3},
  pages={295},
  year={1967},
  publisher={IOP Publishing}
}

@article{mei2018mean,
  title={A mean field view of the landscape of two-layer neural networks},
  author={Mei, Song and Montanari, Andrea and Nguyen, Phan-Minh},
  journal={Proceedings of the National Academy of Sciences},
  volume={115},
  number={33},
  pages={E7665--E7671},
  year={2018},
  publisher={National Academy of Sciences}
}

@article{chizat2018global,
  title={On the global convergence of gradient descent for over-parameterized models using optimal transport},
  author={Chizat, Lenaic and Bach, Francis},
  journal={Advances in neural information processing systems},
  volume={31},
  year={2018}
}

\newpage

\appendix

\section{Neural Estimation of EIGW}
\label{appen:EIGW}
Similarly to the quadratic case, EIGW decomposes as (cf. Section 2.2 in \cite{rioux2024entropic}) $\mathsf{IGW}^{\eps}=\sF_1+\sF_2^{\eps}$, where 
\begin{align*}
    \sF_1(\mu,\nu)&=\int |\langle x,x'\rangle|^2d\mu\otimes \mu(x,x')+\int |\langle y,y'\rangle|^2d\nu\otimes \nu(y,y'),
    \\
    \sF_2^{\eps}(\mu,\nu)&=\inf_{\pi\in\Pi(\mu,\nu)} -2\int \langle x,x'\rangle\langle y,y'\rangle \,d\pi\otimes \pi(x,y,x',y')+\eps \KL(\pi\|\mu\otimes \nu) ,
\end{align*}
with the distinction that the above decomposition does not require the distributions to be centered. Observe that $\mathsf{IGW}^{\eps}$ is invariant under transformations induced by orthogonal matrices, but not under translations~\cite{le2022entropic}. The same derivation underlying Theorem 1 in \cite{zhang2024gromov} results in a dual form for $\sF_2^{\eps}$. 

\begin{lemma}[EIGW duality; Lemma 2 in \cite{rioux2024entropic}]
\label{lemma:F2Dual}
    Fix $\eps>0$, $(\mu,\nu)\in\cP_2(\RR^{d_x})\times\cP_2(\RR^{d_y})$ and let $M_{\mu,\nu}\coloneqq \sqrt{M_2(\mu)M_2(\nu)}$. Then, 
    \begin{equation}
    \label{eq:F2Decomp}
       \mathsf{IGW}^{\eps}(\mu,\nu)=\sF_1(\mu,\nu)+\inf_{\bA \in\RR^{d_x\times d_y}}\Big\{8\|\bA \|_\F^2+\IOT_\bA ^{\mspace{1mu}\eps}(\mu,\nu)\Big\},
     \end{equation}
    where $\IOT_\bA ^{\mspace{1mu}\eps}$ is the EOT problem with cost $c_{\bA}:(x,y)\in\RR^{d_x}\times \RR^{d_y}\mapsto-8x^{\intercal}\bA  y$. 
    Moreover, the infimum is achieved at some $\bA _{\star}\in\cD_{M_{\mu,\nu}}\coloneqq [-M_{\mu,\nu}/2,M_{\mu,\nu}/2]^{d_x\times d_y}$.
\end{lemma}

As in the quadratic case, we can also establish a regularity theory for semi-dual EOT potential of $\IOT_\bA ^{\mspace{1mu}\eps}(\mu,\nu)$, uniformly in $\bA$. In particular, Lemmas \ref{lemma:regularity_EOT} and \ref{lemma:c_transform_regularity} from the proof of \cref{thm:Bound_EGW1} would go through with minimal modification. This enables the same neural estimation approach for the EIGW cost and alignment plan. 

Analogously to \eqref{eq:NE_EGW}, define the EIGW NE as:
\begin{equation}
\label{eq:NE_EIGW}\mathsf{IGW}^\eps_{k,a}(X^n, Y^n)=\sF_1(X^n, Y^n)+\inf _{\bA  \in \cD_{M }}\left\{8\|\bA \|_\F^2+\widehat\IOT_{\bA , k,a}(X^n, Y^n)\right\},
\end{equation}
where $\widehat\IOT_{\bA, k,a}(X^n, Y^n)$ is NE of EOT with cost function $c_{\bA}(x,y)=-8x^{\intercal}\bA y$, given by
\begin{equation}\label{eq:NE_IOT}
     \widehat\IOT_{\bA , k, a}^{\mspace{1mu}\eps}(X^n, Y^n)\coloneqq 
   \max_{f \in \cF_{k, d_x}(a)} \frac{1}{n}\sum_{i=1}^nf(X_i)-\frac{\eps}{n}\sum_{j=1}^n\log\left(\frac{1}{n}\sum_{i=1}^n\exp\left(\frac{f(X_i)+8X_i^{\intercal}\bA Y_j}{\eps}\right)\right).
\end{equation}

Following the same analysis as in quadratic case, we have the following uniform bound of the effective error in terms of the NN and sample sizes for EIGW distances. The below result in analogous to \cref{thm:Bound_EGW1} from the main text that treats the quadratic EGW distance. 

\begin{theorem}[EIGW cost neural estimation; bound 1] \label{thm:Bound_EIGW}
     There exists a constant $C>0$ depending only on $d_x,d_y$, such that setting $a=C(1+\eps^{1-s})$ with $s=\left\lfloor d_x  / 2\right\rfloor+3$, we have
\begin{equation}
\begin{aligned}
&\sup _{(\mu, \nu) \in \mathcal{P}(\cX) \times \mathcal{P}(\cY)} \mathbb{E}\left[\left|\widehat{\mathsf{IGW}}_{k,a}^{\mspace{1mu}\eps}(X^n, Y^n)-\mathsf{IGW}^{\eps}(\mu, \nu)\right|\right]\\
&\qquad\qquad\lesssim_{d_x,d_y} \left(1+\frac{1}{\eps^{\left\lfloor \frac{d_x}{2}\right\rfloor+2}}\right) k^{-\frac{1}{2}}+\min\left\{1+\frac{1}{\eps^{\left\lceil d_x+\frac{d_y}{2}\right\rceil+4}}\,,\left(1+\frac{1}{\eps^{\left\lfloor \frac{d_x}{2}\right\rfloor+2}}\right)\sqrt{k}\right\} n^{-\frac{1}{2}}.
\end{aligned}
\end{equation}
\end{theorem}
\cref{thm:Bound_EIGW} can be derived via essentially a verbatim repetition of the proof of \cref{thm:Bound_EGW1}. The only difference is that the constant $C$ slightly changes since we are working with the cost function $c_{\mathbf{A}}(x,y)=-8 x^{\intercal} \mathbf{A} y$ (which differs from the one in the quadratic case). Similarly, we can also establish analogous results to Theorems \ref{thm:Bound_EGW2} and \ref{thm:Bound_NeuralCoupling} for the EIGW setting. The statements and their derivations remain exactly the same, and are omitted for brevity.

\section{Proofs of Technical Lemmas}\label{appen:proofs_lemma}
\subsection{Proof of Lemma \ref{lem:delta-oracle}}
\label{appen:proof:delta-oracle}
Define the softmax function as $$\operatorname{softmax}(z):=\left(\frac{e^{z_1}}{\sum_{i=1}^ne^{z_i}},\cdots, \frac{e^{z_n}}{\sum_{i=1}^ne^{z_i}}\right)$$ for $z\in\mathbb{R}^n$. We first prove that the softmax function satisfies:
    \[
    \|\operatorname{softmax}(z+\Delta)-\operatorname{softmax}(z)\|_{\infty} \leq \frac{1}{2}\|\Delta\|_{\infty} .
    \]
    for any $z, \Delta \in \mathbb{R}^n$. The Jacobian of softmax function is:
    \[
    J(z)=\nabla_z \operatorname{softmax}(z)=\operatorname{diag}(p)-p p^{\top},
    \]
    where $p=\operatorname{softmax}(z)$. Then
    \begin{align*}
         \|J(z)\|_{\infty \rightarrow \infty}&:=\sup_{y\neq 0} \frac{\|J(z) y\|_\infty}{\|y\|_\infty}\\
         &=\max _i \sum_{k=1}^n\left|J_{i k}\right|\\
         &=\max _i\left(p_i\left(1-p_i\right)+\sum_{k \neq i} p_i p_k\right)\\
         &=\max _i 2 p_i\left(1-p_i\right)\\
         &\leq \frac{1}{2},
    \end{align*} thus by Mean Value Inequality, we have
    \[
    \|\operatorname{softmax}(z+\Delta)-\operatorname{softmax}(z)\|_{\infty} \leq\left(\sup _{\tilde{z}}\|J(\tilde{z})\|_{\infty \rightarrow \infty}\right)\|\Delta\|_{\infty} \leq \frac{1}{2}\|\Delta\|_{\infty} .
    \]
    Now, recall that our neural coupling is defined as 
    \[\big[\widetilde{\boldsymbol{\Pi}}^\varepsilon_{\mathbf{A}}\big]_{ij}=\frac{\exp \left(\frac{\hat{f}_\mathbf{A} \left(X_i\right)-c_{\mathbf{A}}\left(X_i, Y_j\right)}{\varepsilon}\right)}{n \sum_{k=1}^n \exp \left(\frac{\hat{f}_\mathbf{A} \left(X_k\right)-c_{\mathbf{A}}\left(X_k, Y_j\right)}{\varepsilon}\right)}, \quad i,j=1,\ldots,n.   
    \] 
    Then
    \[
    \widetilde{\boldsymbol{\Pi}}^\varepsilon_{\mathbf{A}}=\frac{1}{n}\left(\operatorname{softmax}\left(z^{(1)}(\hat{f}_{\mathbf{A}})\right),\cdots,\operatorname{softmax}\left(z^{(n)}(\hat{f}_{\mathbf{A}})\right)\right),
    \]
    where $z^{(j)}(\hat{f}_{\mathbf{A}}):=(z_1^{(j)}(\hat{f}_{\mathbf{A}}),\cdots,\cdots,z_n^{(j)}(\hat{f}_{\mathbf{A}}))$ and $z_i^{(j)}(\hat{f}_{\mathbf{A}}):=\frac{\hat{f}_{\mathbf{A}}\left(X_i\right)-c_{\mathbf{A}}\left(X_i, Y_j\right)}{\varepsilon}$. Therefore,
    \begin{align*}
    \left\|\widetilde{\boldsymbol{\Pi}}^\varepsilon_{\mathbf{A}}-\boldsymbol{\Pi}^\varepsilon_{\mathbf{A}}\right\|_{\infty} &=\frac{1}{n} \max_j \left\|\operatorname{softmax}(z^{(j)}(\hat{f}_{\mathbf{A}}))-\operatorname{softmax}(z^{(j)}(\varphi_{\mathbf{A}}))\right\|_{\infty}\\
    &\leq \frac{1}{2n}\|z^{(j)}(\hat{f}_{\mathbf{A}})-z^{(j)}(\varphi_{\mathbf{A}})\|_\infty\\
    & =\frac{\delta^{\prime}}{2 n \varepsilon}:=\delta.
    \end{align*} 
    \qed

\subsection{Proof of Lemma \ref{lemma:regularity_EOT}}
\label{appen:proof:lemma:regularity_EOT}
The existence of optimal potentials follows by standard EOT arguments \cite[Lemma 1]{goldfeld2022limit}. Recall that EOT potentials are unique up to additive constants (see section \ref{subsec:EOT}). Thus, let $\left(\varphi^0, \psi^0\right) \in L^1(\mu) \times L^1(\nu)$ be optimal EOT potentials for the cost $c$, solving dual formulation \eqref{eq:EOT_dual}, and we can assume without loss of generality that $\int \varphi^0 d\mu=\int \psi^0 d\nu=\frac{1}{2} \OT_c^{\mspace{1mu}\eps}(\mu, \nu)$.

    Recall that the optimal potentials satisfies the Schr\"odinger system from \eqref{eq:Schrodinger system}. Define new functions $\varphi$ and $\psi$ as
$$
\begin{array}{ll}
\varphi(x)\coloneqq-\eps\log \int_{\cY} \exp\left(\frac{\psi^0(y)-c(x, y)}{\eps}\right) d\nu(y), & x \in \cX \\
\psi(y)\coloneqq \varphi^{c ,\eps}(y), & y \in \cY.
\end{array}
$$
These integrals are clearly well-defined as the integrands are everywhere positive on $\cX$ and $\cY$, and $\varphi^0, \psi^0$ are defined on the supports of $\mu, \nu$ respectively. Now We show that $\varphi,\psi$ are pointwise finite. For the upper bound, by Jensen's inequality, we have
$$
\varphi(x) \leq \int_{\cY} c(x, y)-\psi^0(y) d \nu(y) \leq \frac{3}{2}\|c\|_{\infty, \mathcal{X} \times \mathcal{Y}},
$$
the second inequality follows from $\int \psi^0 d\nu=\frac{1}{2} \OT_{c}^{\mspace{1mu}\eps}(\mu, \nu)$ and
$$
-\OT_{c}^{\mspace{1mu}\eps}(\mu, \nu) \leq \int_{\cX\times\cY}-c(x,y)
d\gamma_{\star}^\eps \leq \|c\|_{\infty, \mathcal{X} \times \mathcal{Y}},
$$
where $\gamma^\eps_{\star}$ is optimal coupling for $\OT_{c}^{\mspace{1mu}\eps}(\mu, \nu)$. The upper bound holds similarly for $\psi$ on $\cY$ and $\psi^0$ on the support of $\nu$. For lower bound, we have
\[
-\varphi(x) \leq \eps\log\int_{\cY} \exp\left(\frac{\frac{5}{2}\|c\|_{\infty, \mathcal{X} \times \mathcal{Y}}}{\eps}\right) d\nu(y)=\frac{5}{2}\|c\|_{\infty, \mathcal{X} \times \mathcal{Y}}
\]
Note that $\varphi_\bA$ is defined on $\cX$, with pointwise bounds proven above. By Jensen's inquality,
$$
\begin{aligned}
&\int_{\cX}\left(\varphi^0-\varphi\right) d \mu+\int_{\cY}\left(\psi^0-\psi\right) d\nu \\
& \leq \eps\log \int_{\cX} \exp\left(\frac{\varphi^0-\varphi}{\eps}\right) d\mu+\log \int_{\cY} \exp\left(\frac{\psi^0-\psi}{\eps}\right) d \nu \\
& =\eps\log \int_{\cX \times \cY} \exp\left(\frac{\varphi^0(x)+\psi^0(y)-c(x, y)}{\eps}\right) d \mu \otimes \nu \\
& \qquad\qquad\qquad\qquad\qquad\qquad+\eps\log \int_{\cX \times \cY} \exp\left(\frac{\varphi(x)+\psi^0(y)-c(x, y)}{\eps}\right) d \mu \otimes \nu\\
& =0 .
\end{aligned}
$$
Since $\left(\varphi^0, \psi^0\right)$ maximizes \eqref{eq:EOT_dual}, so does $(\varphi, \psi)$ and thus they are also optimal potentials. Therefore, $\varphi$ solves semi-dual formulation \eqref{eq:EOT_semidual}. By the strict concavity of the logarithm function we further conclude that $\varphi=\varphi^0$ $ \mu$-a.s and $\psi=\psi^0$ $\nu$-a.s.

The differentiability of $(\varphi, \psi)$ is clear from their definition. For any multi-index $\alpha$, the multivariate Faa di Bruno formula (see \cite[Corollary 2.10]{constantine1996multivariate}) implies
\begin{equation}
\label{eq:constant_C_s_1}
    -D^\alpha \varphi(x)=\eps\sum_{r=1}^{|\alpha|} \sum_{p(\alpha, r)} \frac{\alpha !(r-1) !(-1)^{r-1}}{\prod_{j=1}^{|\alpha|}\left(k_{j} !\right)\left(\beta_{j} !\right)^{k_j}} \prod_{j=1}^{|\alpha|}\left(\frac{D^{\beta_j} \int_{\cY} \exp(\frac{\psi^0(y)-c(x, y)}{\eps}) d \nu(y)}{\int_{\cY} \exp(\frac{\psi^0(y)-c(x, y)}{\eps}) d \nu(y)}\right)^{k_j},
\end{equation}
where $p(\alpha, r)$ is the collection of all tuples $\left(k_1, \cdots, k_{|\alpha|} ; \beta_1, \cdots, \beta_{|\alpha|}\right) \in \mathbb{N}^{|\alpha|} \times \mathbb{N}^{d_x \times|\alpha|}$ satisfying $\sum_{i=1}^{|\alpha|} k_i=r, \sum_{i=1}^{|\alpha|} k_i \beta_i=\alpha$, and for which there exists $s \in\{1, \ldots,|\alpha|\}$ such that $k_i=0$ and $\beta_i=0$ for all $i=1, \ldots|\alpha|-s, k_i>0$ for all $i=|\alpha|-s+1, \ldots,|\alpha|$, and $0 \prec \beta_{|\alpha|-s+1} \prec \cdots \prec \beta_{|\alpha|}$. For a detailed discussion of this set including the linear order $\prec$, please refer to \cite{constantine1996multivariate}. For the current proof we only use the fact that the number of elements in this set solely depends on $|\alpha|$ and $r$. Given the above, it clearly suffices to bound $\left|D^{\beta_j} \int_{\cY} \exp\left(\psi^0(y)-c(x, y)\right) d\nu(y)\right|$. First, we apply the same formula to $D^{\beta_j} e^{-c(x, y)/\eps}$ and obtain
\begin{equation}
\label{eq:constant_C_s_2}
   D^{\beta_j} e^{\frac{-c(x, y)}{\eps}}=\sum_{r^{\prime}=1}^{\left|\beta_j\right|} \sum_{p\left(\beta_j, r^{\prime}\right)} \frac{(-1)^{-r^\prime}\beta_{j} ! \eps^{-r^\prime}}{\prod_{i=1}^{\left|\beta_j\right|}\left(k_j^{\prime} !\right)\left(\eta_{i} !\right)^{k_i^{\prime}}} e^{\frac{-c(x, y)}{\eps}} \prod_{i=1}^{\left|\beta_j\right|}\left(D^{\eta_i}c(x,y)\right)^{k_i^{\prime}} ,
\end{equation}
where $p\left(\beta_j, r^{\prime}\right)$ is a set of tuples $\left(k_1^{\prime}, \ldots, k_{\left|\beta_j\right|}^{\prime} ; \eta_1, \ldots, \eta_{\left|\beta_j\right|}\right) \in \mathbb{N}^{\left|\beta_j\right|} \times \mathbb{N}^{d_x \times\left|\beta_j\right|}$ defined similarly to the above. Since $c(x,y)$ is smooth, we can conclude that $\left|D^{\eta_i} c(x, y)\right|\leq C_{c,\|\cX\|,\|\cY\|}$ for some constant $C_{c,\|\cX\|,\|\cY\|}$ that depends on $c(x,y)$ and the diameter of $\cX$ and $\cY$. 
Consequently, for $0<\eps< 1$, we have
$$
    \left|\frac{D^{\beta_j} \int_{\cY} \exp(\frac{\psi^0(y)-c(x, y)}{\eps}) d \nu(y)}{\int_{\cY} \exp(\frac{\psi^0(y)-c(x, y)}{\eps}) d \nu(y)}\right| \leq C_{\beta_j}\eps^{-|\beta_j|}C^{|\beta_j|}_{c,\|\cX\|,\|\cY\|},
$$
and for $\eps \geq1$,
$$
    \left|\frac{D^{\beta_j} \int_{\cY} \exp(\frac{\psi^0(y)-c(x, y)}{\eps}) d \nu(y)}{\int_{\cY} \exp(\frac{\psi^0(y)-c(x, y)}{\eps}) d \nu(y)}\right| \leq C_{\beta_j}\eps^{-1}C^{|\beta_j|}_{c,\|\cX\|,\|\cY\|}.
$$
Plugging back, we obtain 
$$
\left|D^\alpha \varphi(x)\right| \leq C_{|\alpha|}(1+\eps^{1-|\alpha|})C^{|\alpha|}_{c,\|\cX\|,\|\cY\|}.
$$
Analogous bound holds for $\psi$.

In particular, when we take $c=c_\bA$ for any $\bA\in\mathcal{D}_M$, there exist semi-dual EOT potentials $(\varphi_{\bA},\varphi^{c,\eps}_{\bA})$ for $\OT_{\bA}^{\mspace{1mu}\eps}(\mu, \nu)$, such that
\begin{equation}
\label{eq:regularity_EGW_potential}
\begin{aligned}
    \|\varphi_{\bA}\|_{\infty,\cX} &\leq 6 d_x d_y+40 (d_xd_y)^{3/2}\\
    \left\|D^\alpha \varphi_{\bA}\right\|_{\infty,\cX} &\leq C_{s}\left(1+\eps^{1-s}\right)\left(1+d_y\right)^{s}\left(1+\sqrt{d_x} d_y+\sqrt{d_x}\right)^{s},\quad \forall \alpha \in \mathbb{N}_0^{d_x} \text { with } 1 \leq|\alpha| \leq s.
\end{aligned}
\end{equation}
For the first bound, we use the fact $\|c_\bA\|_{\infty, \mathcal{X} \times \mathcal{Y}}\leq4\|x\|^2\|y\|^2+32\|x\|\|A\|_{\F}\|y\|\leq4 d_x d_y+16 M d_x d_y$. For derivatives bound, observe that
\begin{align*}
   \left|D^{\eta_i} c_{\bA}(x, y)\right|=\left|D^{\eta_i}\left(4\|x\|^2\|y\|^2+32 x^{\intercal} \bA  y\right)\right| &\leq 8\left(1+\|y\|^2\right)\left(1+\sqrt{d_y} M+\|x\|\right) \\
   &\leq 8(1+d_y)(1+\sqrt{d_y}M +\sqrt{d_x}),
\end{align*}
where the first inequality follows from proof of Lemma 3 in \cite{zhang2024gromov}. We obtain 
$$
\left|D^\alpha \varphi_\bA(x)\right| \leq C_{|\alpha|}(1+\eps^{1-|\alpha|})(1+d_y)^{|\alpha|}(1+\sqrt{d_y}M +\sqrt{d_x})^{|\alpha|}.
$$
The proof is completed by plugging the definition of $M $ into these bounds.
\qed

\subsection{Proof of Lemma \ref{lemma:neural_approx_reduction}}
\label{appen:proof_neural_approx_reduction}
For any $f\in \cF_{k,a}(a)$, we know that $\|f\|_{\infty,\cX} \leq 3 a(\|\cX\|+1)=6a$, so NNs are uniformly bounded. This implies that $\OT_c^{\eps}(\mu,\nu) \geq \OT^{\eps}_{k,a}(\mu,\nu)$. Since $\varphi$ satisfies \eqref{eq:regularity_EOT}, it's uniformly bounded on $\cX$. Then, the following holds:
\begin{align*}
\left|\OT_c^{\eps}(\mu,\nu)-\OT^{\eps}_{c,k,a}(\mu,\nu)\right|  =&\OT_c^{\eps}(\mu,\nu)-\OT^{\eps}_{c,k,a}(\mu,\nu) \\
\leq &\mathbb{E}_{\mu} |\varphi-f|+\mathbb{E}_{\nu} |\varphi^{c,\eps}-f^{c,\eps}|\\
\leq & 2\|\varphi-f\|_{\infty,\cX}.
\end{align*}  
The last inequality holds by an observation that, 
\begin{equation}
\label{eq:lipshcitz_c_tranform}
    |\varphi^{c,\eps}(y)-f^{c,\eps}(y)| \leq \left\|\varphi-f\right\|_{\infty, \cX},\quad \forall y\in \cY.
\end{equation}
Indeed, note that for any $y \in \cY$,
\begin{align*}
    &\varphi^{c,\eps}(y)-f^{c,\eps}(y)\\
    &= -\eps \log \left(\frac{\int \exp \left(\frac{\varphi(x)-c(x, y)}{\eps}\right) d \mu(x)}{\int \exp \left(\frac{f(x)-c(x, y)}{\eps}\right) d \mu(x)}\right)\\
    & = -\eps \log \left(\frac{\int \exp\left(\frac{\varphi(x)-f(x)}{\eps}\right)\exp \left(\frac{f(x)-c(x, y)}{\eps}\right) d \mu(x)}{\int\exp \left(\frac{f(x)-c(x, y)}{\eps}\right) d \mu(x)}\right)\\
    & \geq -\eps \log \left(\frac{\int \exp\left(\frac{\|\varphi-f\|_{\infty,\cX}}{\eps}\right)\exp \left(\frac{f(x)-c(x, y)}{\eps}\right) d \mu(x)}{\int\exp \left(\frac{f(x)-c(x, y)}{\eps}\right) d \mu(x)}\right)\\
    & = -\left\|\varphi-f\right\|_{\infty, \cX}.
\end{align*}
Similarly, we can have that
\begin{align*}
    &\varphi^{c,\eps}(y)-f^{c,\eps}(y) \\
    & \leq -\eps \log \left(\frac{\int\exp\left(\frac{-\|\varphi-f\|_{\infty,\cX}}{\eps}\right)\exp \left(\frac{f(x)-c(x, y)}{\eps}\right) d \mu(x)}{\int \exp \left(\frac{f(x)-c(x, y)}{\eps}\right) d \mu(x)}\right)\\
    & = \left\|\varphi-f\right\|_{\infty, \cX}.
\end{align*}
\qed

\end{document}